\documentclass[11pt,a4paper,reqno]{amsart}
\usepackage{verbatim}
\usepackage{graphicx}
\usepackage[utf8]{inputenc}
\usepackage{amsmath,amsfonts,amssymb,amscd,amsthm,txfonts,hyperref}
\usepackage{xcolor}
\newtheorem{theorem}{Theorem}[section]
\newtheorem{proposition}[theorem]{Proposition}
\newtheorem{lemma}[theorem]{Lemma}
\newtheorem{corollary}[theorem]{Corollary}
\newtheorem{remark}{Remark}[section]

\theoremstyle{definition}
\newtheorem{definition}[theorem]{Definition}
\newtheorem{notation}[theorem]{Notation}

\newcommand{\R}{\mathbb{R}}
\newcommand{\C}{\mathbb{C}}
\newcommand{\N}{\mathbb{N}}
\newcommand{\Z}{\mathbb{Z}}

\newcommand{\lap}{\bigtriangleup}

\newcommand{\E}{\mathbb E}
\renewcommand{\Re}{\mbox{Re}}

\newcommand{\an}[1]{\langle #1 \rangle}
\newcommand{\grad}{\bigtriangledown}
\newcommand{\mat}[2]{\begin{pmatrix} #1 \\ #2 \end{pmatrix}}
\newcommand{\be}{\begin{equation}}
\newcommand{\ee}{\end{equation}}
\newcommand{\bee}{\begin{eqnarray*}}
\newcommand{\eee}{\end{eqnarray*}}
\newcommand{\lab}{\label}
\newcommand{\fref}{\eqref}
\newcommand{\ba}{\begin{array}}
\newcommand{\ea}{\end{array}}
\newcommand{\pa}{\partial}
\newcommand{\e}{\varepsilon}
\newcommand{\bea}{\begin{eqnarray}}
\newcommand{\eea}{\end{eqnarray}}
\newcommand{\non}{\nonumber}

\begin{document}

\title{Stability of equilibria for a Hartree equation \\ for random fields}


\author[C. Collot]{Charles Collot}
\address{Courant Institute of Mathematical Sciences, New York University, New York, United States of America}
\email{cc5786@nyu.edu}
\author[A.-S. de Suzzoni]{Anne-Sophie de Suzzoni}
\address{CMLS, \'Ecole polytechnique, CNRS, Universit\'e Paris-Saclay, 91128 Palaiseau Cedex, France.}
\email{anne-sophie.de-suzzoni@polytechnique.edu}

\begin{abstract} 
We consider a Hartree equation for a random variable, which describes the temporal evolution of infinitely many Fermions. On the Euclidean space, this equation possesses equilibria which are not localised. We show their stability through a scattering result, with respect to localised perturbations in the defocusing case in high dimensions $d\geq 4$. This provides an analogue of the results of Lewin and Sabin \cite{LS2}, and of Chen, Hong and Pavlovi\'c \cite{CHP2} for the Hartree equation on operators. The proof relies on dispersive techniques used for the study of scattering for the nonlinear Schr\"odinger and Gross-Pitaevskii equations.
\end{abstract}

\maketitle


\section{Introduction}

\subsection{Mean-field dynamics of an infinite number of Fermions}

The present work concerns the following Hartree equation for a random variable:
\begin{equation}\label{eqonrv}
i\partial_t X = -\Delta X +\left( w* \E(|X|^2)\right) X.
\end{equation}
Here $X:I\times \mathbb R^d\times \Omega \rightarrow \C $ is a time-dependent random field over the Euclidean space $\R^d$, and $(\Omega,\mathcal A,d\omega)$ is the underlying probability space. The expectation is with respect to this probability space: $\E (|X|^2)(t,x):=\int_\Omega |X|^2(t,x,\omega)d\omega$. The convolution product is denoted by $*$ and $w$ is a real-valued pair interaction potential. By this, we mean that we only consider interactions between two particles (and no more), and that this interaction is characterized by $w$. Equation \fref{eqonrv} has been introduced in \cite{dS} as an effective dynamics for a large, possibly infinite, number of Fermions in a mean field regime.\\

\noindent Indeed, consider the evolution of a finite number of Fermions interacting through the potential $w$. Under some mean-field hypothesis, as the number of particles tend to infinity, the system is approximated to leading order by a system of $N$ coupled Hartree equations on $\mathbb R \times \mathbb R^d$:
\be \lab{eq:hartreefinite}
i\pa_t u_j =-\Delta u_j+\left(w* (\sum_{k=1}^N |u_k|^2)\right)u_j, \ \ j=1,...,N,
\ee
where the family $(u_j)_{1\leq j\leq N}$ is required to be orthonormal to be compliant with the Pauli principle (which is preserved by the dynamics). Let us associate to this orthonormal family the operator $\gamma =\sum_{1}^N |u_j \rangle \langle u_j|$, which is the orthogonal projection onto $\text{Span}\{(u_j)_{1\leq j \leq N} \}$. There exists a large literature about the derivation of this system of equations and about other related approximation results. In particular, it has been showed that, if the wave function of the original fermionic system is close to a Slater determinant, then, in the mean-field limit and under sufficient conditions for $w$, the associated one-particle density matrix converges to the above operator $\gamma$. For the derivation of Equation \ref{eq:hartreefinite} from many body quantum mechanics we refer to \cite{FK,EESY,BPS,BGGM,
BEGMY,BJPSS}. Note that the so-called exchange term appearing in the Hartree-Fock equation is not present in \fref{eq:hartreefinite}, which is motivated by the fact that it is of lower order in certain regimes, see the aforementioned references.\\

\noindent To deal with infinitely many particles, it is customary to use the density matrices framework, which is an operator formalism. Namely, the family $(u_j)_{1\leq j \leq N}$ solves Equation \fref{eq:hartreefinite} if and only if the operator $\gamma$ defined above solves the corresponding Hartree equation:
\be \lab{eq:hartreeoperator}
i\pa_t \gamma = [-\Delta +w*\rho_{\gamma},\gamma].
\ee
Above, $[,]$ denotes the commutator, and $\rho_{\gamma}(x)=\tilde \gamma (x,x)$ is the density of particles, that is the diagonal of the integral kernel $\tilde \gamma (x,y) $ of the operator $\gamma$. An infinite number of particles can then be modelled by a solution of \fref{eq:hartreeoperator} which is not of finite rank (the rank of the operator being, by the derivation of the model, the number of particles). Solutions of \fref{eq:hartreeoperator} with an infinite number of particles were studied previously in \cite{BPF,BPF2,C,Z} for exemple, and more recently in \cite{LS,LS2,CHP,CHP2}.\\

\noindent In \cite{dS}, the second author proposed \fref{eqonrv} as an alternative equation to \fref{eq:hartreeoperator}. It generalises Equation \fref{eq:hartreefinite} for a finite number of particles in the following sense. To an orthonormal family $(u_j)_{1\leq j \leq N}$, one can associate the random variable $X(x,\omega)= N^{-1/2} \sum_1^N u_j(x)g_j(\omega)$, where $(g_j)_{1\leq j\leq N}$ is an orthonormal family in $L^2_\omega$. The family $(u_j)_{1\leq j \leq N}$ then solves \fref{eq:hartreefinite} if and only if the random variable $X$ solves \fref{eqonrv}. Equation \fref{eqonrv} is also in close correspondance with the Hartree equation for density matrices \fref{eq:hartreeoperator}. We refer to \cite{dS} for how to relate the solutions of the two Equations \fref{eqonrv} and \fref{eq:hartreeoperator} and their corresponding equilibria. One reason behind the study of \fref{eqonrv} is that this equation shares more direct resemblances with the commonly studied nonlinear Schr\"odinger equation.\\

\noindent The existence of solutions of Equation \fref{eqonrv} in $L^2_\omega H^s_x$ is investigated in \cite{dS}. The local well-posedness is established in the case of localised initial data, as well as in the case of localised perturbations of the equilibria described below. In particular, almost everywhere in the probability space, the random variable solves the corresponding Schr\"odinger equation in integrated formulation. Though the results are stated and proved in the case of a Dirac potential $w=\delta_{\{x=0 \}}$, their adaptation to the present case of regular interaction potentials is straightforward. Moreover, as for the defocusing nonlinear Schr\"odinger equation, scattering for large but localised solutions is expected for Equation \fref{eqonrv}, at least in energy critical and subcritical regimes. This has been showed in dimension 3 in the case of a Dirac potential in \cite{dS}. Nontrivial equilibria are thus non-localised, which corresponds to being not of trace-class in the framework of density matrices.\\


\subsection{Statement of the result}

\noindent The equilibria at stake in the present paper are the following. For Equation \fref{eq:hartreeoperator}, any non-negative Fourier multiplier $\gamma:u\mapsto \mathcal F^{-1}(|f(\xi)|^2\mathcal F u)$ with symbol $|f|^2$ is a stationary solution, where $\mathcal F$ denotes the Fourier transform. Note that density operators are non-negative, we write the symbol in the form $|f|^2$ to be able to give an analogous equilibrium in the random framework. The analogous equilibrium for Equation \fref{eqonrv} are given by Wiener integrals
\be \lab{id:Y}
Y_f(t,x,\omega):= \int_{\xi \in \mathbb R^d} f(\xi) e^{i \xi.x-it(m+|\xi|^2)}dW(\xi),
\ee
for a distribution function $f:\mathbb R^d \rightarrow \mathbb C$ (note that this equilibrium has the same law if we replace $f$ by $|f|$ so that we can assume $f:\mathbb R^d\rightarrow [0,+\infty)$). Above $dW(\xi)$ denotes infinitesimal complex Gaussians characterised by
$$
\E \left(\overline{dW(\eta)}dW(\xi)\right)=\delta_{\eta-\xi}d\eta d\xi,
$$
and the scalar $m$ is given by $m:=\int_{\mathbb R^d} w(x)dx \int_{\mathbb R^d} f(\xi)d\xi$. We refer to \cite{S} for more information an random Gaussian fields. The function $Y_f$ is a solution of Equation \fref{eqonrv}. It is not stationary but almost, since its law is invariant by time and space translations (and in particular is not localised).\\

\noindent In the seminal work \cite{LS2}, the authors show the stability of the above equilibria for the Equation \fref{eq:hartreeoperator} for density matrices in dimension 2. Important tools are dispersion estimates for orthonormal systems \cite{FLLS,FS}. This work has been extended to higher dimension in \cite{CHP2}. Note that in higher dimension, some structural hypothesis is made on the interaction potential $w$, to solve some technical difficulties about a singularity in low frequencies of the equation that we will identify precisely in the sequel. The stability result corresponds to a scattering property in the vicinity of these equilibria: any small and localised perturbation evolves asymptotically into a linear wave which disperses. We mention equally \cite{LS,CHP} about problems of global well-posedness for the equation on density matrices.\\

\noindent The problem of the stability of the equilibria \fref{id:Y} for Equation \fref{eqonrv} shares similarities with the stability of the trivial solution for the Gross-Pitaevskii equation $i\pa_t \psi=-\Delta \psi+(|\psi|^2-1)\psi$. In both problems the linearised dynamics has distinct dispersive properties at low and high frequencies, making the nonlinear stability problem harder, especially in low dimensions. The proof of scattering for small data for the Gross Pitaevskii equation was done in \cite{GNT,GNT2,GNT3,GHN}. We here use spaces with different regularities at low and high frequencies, inspired by \cite{GNT}. \\

\noindent The result of the present work is the stability of the equilibria \fref{id:Y} for Equation \fref{eqonrv}, via the proof of scattering of perturbations to linear waves in their vicinity. Our techniques however differs from those used in \cite{LS2,CHP2}, see the strategy of the proof after Theorem \ref{thmmain}. We hope that the present proof provides different insights than the ones in the framework of operators, as well as relaxing some of the hypotheses on the potential $w$.\\

\noindent In what follows, $\langle \xi \rangle =(1+|\xi|^2)^{1/2}$ denotes the usual Japanese bracket, and we write with an abuse of notation $f(\xi)=f(r)$ with $r=|\xi|$, if $f$ has spherical symmetry. The space $L^2_\omega H^{s}$ is the set of measurable functions $Z:\mathbb R^d\times \Omega \rightarrow \mathbb C$ such that $Z(\cdot,\omega)\in H^s$ almost surely and 
$$
\int_{\mathbb R^d\times \Omega} \langle \xi \rangle^{2s} |\hat Z(\xi,\omega)|^2d\xi d \omega<+\infty.
$$

\begin{theorem}[Stability of equilibria for Equation \fref{eqonrv}]
\label{thmmain}
Let $d\geq 4$. Let $s = \frac{d}2 - 1$. Let $f$ be a spherically symmetric function in $L^2\cap L^\infty(\R^d)$. 
Assume : \begin{itemize} 
\item $\an{\xi}^{\lceil s \rceil} f \in L^2 (\R^d)$,
\item $\int_{\R^d} |\xi|^{1-d}|f\grad  f| < \infty$,
\item $\partial_r |f|^2 <0$ for $r>0$,
\item writing $h$ the inverse Fourier transform of $|f|^2$, $\an{\xi}^2 \pa^{\alpha} h \in L^\infty(\R^d)$ for all $\alpha \in \mathbb N^d$ with $|\alpha|\leq 2 \lceil s \rceil$,
\item $|\xi|^{1-d}(h + \grad h) \in L^1(\R^d)$.
\end{itemize}
There exists $C(f)$ such that for all $w \in W^{s,1}$ that satisfies
$$
\an \xi \hat w \in L^{2(d+2)/(d-2)}\textrm{ and }\|(\hat w)_-\|_{L^\infty} + \hat{w}(0)_+ \leq C(f),
$$
there exists $\epsilon = \epsilon(f,w)$ such that for any $Z_0\in L^2_\omega H^{d/2-1}\cap L^{2d/(d+2)}_xL^2_\omega$ with
$$
\| Z_0 \|_{L^2_\omega H^{d/2-1}}+\| Z_0 \|_{L^{2d/(d+2)}_xL^2_\omega}\leq \epsilon,
$$
the solution of \fref{eqonrv} with initial datum $X_0 = Y_f(t=0) +Z_0$ is global. Moreover, it scatters to a linear solution in the sense that there exists $Z_-,Z_+\in L^2_\omega H^{d/2-1}$ such that
$$
X(t)= Y_f(t) + e^{i(\Delta-m)t}Z_{\pm}+o_{L^2_\omega H^{d/2-1}}(1) \ \ \text{as} \ t \rightarrow \pm \infty.
$$

\end{theorem}

\vspace{0.3cm}

\begin{remark} The conditions on $f$ are satisfied by thermodynamical equilibria for bosonic or fermionic gases at a positive temperature $T$ :
$$
|f(\xi)|^2=\frac{1}{e^{\frac{|\xi|^2-\mu}{T}}-1}, \ \ \mu <0 \ \ \text{ and } \ \ |f(\xi)|^2=\frac{1}{e^{\frac{|\xi|^2-\mu}{T}}+1}, \ \ \mu\in \mathbb R,
$$
respectively, but it is not the case of the fermionic gases of zero temperature : 
$$
|f(\xi)|^2={\bf 1}_{|\xi|^2\leq \mu}, \ \ \mu>0.
$$

\end{remark}

\begin{remark} 

The smallness assumption on $\|(\hat w)_-\|_{L^\infty} $, corresponds to the fact that the equation is not too focusing. The one on $(\hat w(0))_+$ enables the equation \eqref{eqonrv} linearized around $Y_f$ to have enough dispersion. Note that these assumptions appear both in \cite{LS2} and \cite{CHP2}.

\end{remark}

\begin{remark} 

The stability is not a consequence of a strong convergence to $0$ of the perturbation in $H^{d/2-1}$ almost everywhere in probability, but of the dispersion for the linear dynamics implying in particular local convergence. The solution converges back to the same equilibrium $Y_f$. In particular, it does not trigger modulationnal instabilities.

\end{remark}

\begin{remark} In dimension $2$ and $3$, the treatment of the problem has to be somewhat modified compared to the proof we present. The reason for this is some singularity in the low frequencies that can be dealt with thanks to dispersion estimates on dimension higher than $4$, but which are more problematic in dimension $2$ and $3$. Nevertheless, this singularity is more an artefact of the proof than a property of the equation. It can be dealt with by iterating the wave operator twice, computing exactly different terms of the equation to cancel out the singularity. This is related to the fact that in the context of operators, the wave operator has to be iterated more than once to be able to apply dispersion estimates (we refer to \cite{LS2}). In dimension $3$ and in our context, the result would be a bit weaker than in dimension higher than $4$ : we would have to start form an initial datum of typical Sobolev regularity $\frac12$ and get the scattering in $L^2$. \end{remark}

\begin{remark} Even in the defocusing case, some other equilibria can have instabilities. Two plane waves, which are orthogonal in probability, propagating in opposite directions, are linearly instable, which is showed In Section \ref{sec:2w}. Note that the equilibria of Theorem \ref{thmmain} can be seen as a superposition of infinitely many plane waves propagating in different directions, hence this shows the importance of regularity of the underlying function $f$. \end{remark}

\noindent The strategy of the proof is the following. First note that the dynamics of Equation \fref{eqonrv} near the above equilibria is somewhat similar to that of the Gross-Pitaevskii equation. We use a more direct fixed point argument which does not involve iterations of the wave operator as in \cite{LS2,CHP2}, and dispersion properties at low and high frequencies which are inspired from \cite{GNT}. \\

\noindent We start by reducing the proof to finding a correct functionnal framework for our contraction argument. Namely, instead of solving an equation for the perturbation $Z = X-Y_f$, we solve a fix point for the perturbation and the induced potential $Z,V$ where $V = \E(|Z|^2) + 2\Re \E(\bar Y_f Z))$, at the same time. The idea behind this is that even if $V$ contains a linear term in $Z$, it behaves more like a quadratic term in $Z$ (in the sense of the Lebesgue spaces to which they belong) and thus, we can put it in better spaces regarding dispersion.\\

\noindent The fixed point is solved in a classical way by finding the right Banach space $\Theta$ for $Z,V$ and proving suitable estimates on the linear and nonlinear terms. The Lipschitz-continuity of the quadratic part is treated in a classical fashion, in the sense that it requires that $\Theta$ is included in classical Lebesgue, Besov, or Sobolev spaces. \\

\noindent The difficulty comes from the linear term. It can be written 
$$
L = \begin{pmatrix} 0 & L_2 \\ 0 & L_1 \end{pmatrix},
$$
see \fref{id:L}, where $L_1$ corresponds to the analogue of the linear term in \cite{LS2,CHP2}. The invertibility and continuity of $1-L_1$ had been dealt with in both papers and their treatment is more or less sufficient for our argument. Note that this term is the linear response of the equilibrium, related to the so-called Lindhard function \cite{L,GV}.\\

\noindent But $L_2$ is where singularities in low frequencies occurs. To get the continuity of $L_2$, $V$ needs to be in a space that compensates this singularity, namely inhomogeneous Besov spaces, with two levels of regularity, one for the low frequencies and one for the high frequencies. But $V$ contains $\E(|Z|^2)$ and we cannot close the fix point argument for $Z$ in a space that compensates singularities in low frequencies. This is where we use bilinear estimates on inhomogeneous Besov spaces coming from the scattering for Gross-Pitaevskii literature, \cite{GNT}.\\

The paper is organised as follows. In Section \ref{sec-toolbox}, we state a few definitions and known results that we use in the rest of the paper. In Section \ref{sec:principle} we set up the fixed point argument. In Section \ref{sec:embeddings} estimates related to some direct embeddings are given. The linear terms are studied in Section \ref{sec:lin} where explicit formulas, continuity estimates and invertibility conditions are obtained. Estimates for the quadratic terms are proven in Section \ref{sec:quad}, and estimates for the source terms are showed in Section \ref{sec:ini}. Theorem \ref{thmmain} is then proved in Section \ref{sec-proofth}. The last Section \ref{sec:2w} is devoted to an instability result in a defocusing case when the partition function is not smooth.

\subsection*{Notations} 

\begin{notation}[Fourier transform] We define the Fourier transform with the following constants : for $g \in \mathcal S$,
$$
\hat g (\xi) = \mathcal F (g) (\xi) = \int_{\R^d} g(x) e^{-ix\xi}dx,
$$
and the inverse Fourier tranform by
$$
\mathcal F^{-1} (g) (x) = (2\pi)^{-d}\int_{\R^d} g(\xi) e^{ix\xi}d\xi.
$$
\end{notation}

\begin{notation}[Time-space norms] For $p,q\in [1,\infty]$, we denote by $L^p_t,L^q_x = L^p,L^q$ the space
$$
L^p(\R,L^q(\R^d)).
$$

For $p,q\in [1,\infty], s\in \R$, we denote by $L^p_t,W^{s,q}_x = L^p,W^{s,q}$ the space
$$
L^p(\R,W^{s,q}(\R^d)).
$$
In the case $q=2$ we also write it $L^p,H^s$ or $L^p_t,H^s_x$.

When $p=q$, we may write $L^p_{t,x}$ for $L^p_t,L^p_x$.

For $p,q\in [1,\infty], s,t\in \R$, we denote by $L^p_t,(B_q^{s,t})_x = L^p,B_q^{s,t}$ the space
$$
L^p(\R,B_q^{s,t}(\R^d)).
$$
The proper definition of inhomogeneous Besov spaces is given in Section \ref{sec-toolbox}, Definition \ref{def-inhomBesov}.
\end{notation}

\begin{notation}[Probability-time-space norms]
For $p,q\in [1,\infty]$, we denote by $L^2_\omega ,L^p_t,L^q_x =L^2_\omega, L^p,L^q$ the space
$$
L^2(\Omega,L^p(\R,L^q(\R^d))).
$$

For $p,q\in [1,\infty], s\in \R$, we denote by $L^2_\omega,L^p_t,W^{s,q}_x = L^2_\omega,L^p,W^{s,q}$ the space
$$
L^2(\Omega,L^p(\R,W^{s,q}(\R^d))).
$$
In the case $q=2$ we also write it $L^2_\omega,L^p,H^s$ or $L^2_\omega,L^p_t,H^s_x$.

For $p,q\in [1,\infty], s,t\in \R$, we denote by $L^2_\omega,L^p_t,(B_q^{s,t})_x = L^2_\omega,L^p,B_q^{s,t}$ the space
$$
L^2(\Omega,L^p(\R,B_q^{s,t}(\R^d))).
$$
\end{notation}

\begin{notation}[Time-space-probability norms]
For $p,q\in [1,\infty], s\in \R$, we denote by $L^p_t,W^{s,q}_x,L^2_\omega = L^p,W^{s,q},L^2_\omega$ the space
$$
(1-\lap_x)^{s/2}L^p(\R,L^{q}(\R^d,L^2(\Omega))).
$$
In the case $q=2$ we also write it $L^p,H^s,L^2_\omega$, and note that $L^p,H^s,L^2_\omega=L^pL^2_\omega H^s$. In the case $s=0$, we also write it $L^p,L^q,L^2_\omega$.

For $p,q\in [1,\infty], s,t\in \R$, we denote by $L^p_t,(B_q^{s,t})_x,L^2_\omega = L^p,B_q^{s,t},L^2_\omega$ the space induced by the norm
$$
\|g\|_{L^p,B_q^{s,t},L^2_\omega} = \big\|\Big( \sum_{j<0} 2^{2js} \|g_j\|_{L^q_x,L^2_\omega}^2  + \sum_{j\geq 0} 2^{2jt} \|g_j\|_{L^q_x,L^2_\omega}^2 \Big)^{1/2} \big \|_{L^p(\R)}.
$$
\end{notation}

\subsection*{Acknowledgements} 

C. Collot is supported by the ERC-2014-CoG 646650 SingWave. A-S de Suzzoni is supported by ANR 2018 ESSED. Part of this work was done when C. Collot was visiting IH\'ES and he thanks the institute.


\section{Toolbox}\label{sec-toolbox}

In this section, we present existing results in the literature, either classical or more recent ones such as dispersion estimates and matters related to Littlewood-Paley decomposition.

Take $\eta$ a smooth function with support included in the annulus $\{ \xi \in \R^d | \; |\xi| \in (1/2,2)\}$, and define for $j\in \mathbb Z$, $\eta_j(\xi)=\eta(2^{-j}\xi)$. We assume that on $\R^d\smallsetminus \{0\}$, $\sum_j \eta_j = 1$. For any tempered distribution $f \in \mathcal S'(\R)$, we write $f_j = \Delta_j f$ where $\Delta_j$ is the Fourier multiplier by $\eta_j$ i.e. $\hat f_j=\eta_j\hat f$.

\begin{notation}\label{not-LW} We have $f = \sum_j f_j$ and we call this the Littlewood-Paley decomposition of $f$. \end{notation}

\begin{lemma}\label{lem-bernstein}[Bernstein's lemma] Let $a\geq b \geq 1$, there exists $C$ such that for all $j \in \Z$ and all $f\in \mathcal S'$ such that $f_j \in L^b(\R^d)$, we have
$$
\|f_j\|_{L^a} \leq C 2^{jd(\frac1{b}-\frac1{a})} \|f_j\|_{L^b}.
$$
\end{lemma}

We now introduce inhomogeneous Besov spaces where the inhomogeneity comes from a different treatment of high and low frequencies.

\begin{definition}\label{def-inhomBesov}[\cite{GNT}] Let $s,t \in \R$ and $p \geq 1$, we define $B_p^{s,t}$ by the space induced by the norm
$$
\|f\|_{B_p^{s,t}} = \Big( \sum_{j<0} 2^{2js} \|f_j\|_{L^p}^2  + \sum_{j\geq 0} 2^{2jt} \|f_j\|_{L^p}^2 \Big)^{1/2}.
$$
\end{definition}

\begin{remark} This corresponds to taking the homogeneous Besov norm $\dot B_{p,2}^s$ for the low frequencies and $\dot B_{p,2}^t$ for the high frequencies. \end{remark}

We state a few properties of these spaces.

\begin{proposition}\label{prop-Besovembed} Let $s_1 \leq s_2$ and $t_1\geq t_2$ and $p\geq 1$. We have that for all $f \in B_p^{s_1,t_1}$, $f$ also belongs to $B_p^{s_2,t_2}$ and we have
$$
\|f\|_{B_p^{s_2,t_2}} \leq \|f\|_{ B_p^{s_1,t_1}}.
$$
\end{proposition}

\begin{theorem}\label{th-LWth}[Littlewood-Paley theorem] Let $s\geq 0$ and $p \geq 1$. If $p\geq 2$, there exists $C$ such that for all $f \in B_p^{0,s}$ we have 
$$
\|f\|_{W^{s,p}} \leq C \|f\|_{B_p^{0,s}}.
$$
If $p\leq 2$ then there exists $C$ such that for all $f \in W^{s,p}$, we have
$$
\|f\|_{B_p}^{0,s} \leq \|f\|_{W^{s,p}}.
$$
\end{theorem}

We cite here bilinear estimates from \cite{GNT}.

\begin{proposition}\label{prop-bilinBesov}[Lemma 4.1 \cite{GNT}] Let $(s_j)_{1\leq j \leq 3} \in \R^3$, $(t_j)_{1\leq j \leq 3} \in \R^3$, $(p_j)_{1\leq j \leq 3} \in [2,\infty[^3$ such that for all $j \in \{1,2,3\}$,
$$
\max(0,s_j,s_1+s_2+s_3) \leq d\Big( \frac1{p_1} + \frac1{p_2} + \frac1{p_3} -1\Big) \leq t_1+t_2+t_3 \geq t_j  \textrm{ and } s_jp_j<d.
$$
There exists $C$ such that for all $(f_j)_{1\leq j\leq 3}$ such that $f_j \in B_{p_j}^{s_j,t_j}$ for all $j\in \{1,2,3\}$,
$$
\Big| \int_{\R^d} f_1f_2f_3 \Big| \leq C \prod_{j=1}^{3}\|f_j\|_{B_{p_j}^{s_j,t_j}}.
$$
\end{proposition}

We now state some dispersion estimates for the linear flow. Define $S(t) = e^{-it(m-\lap)}$.

\begin{proposition}\label{prop-Striadmiss} Take $p,q \in [2,\infty]$ (note that $d\geq 3$ so that we are excluding the endpoint) such that 
$$
\frac2{p} + \frac{d}{q} = \frac{d}2.
$$
There exists $C$ such that for all $u \in L^2$,
$$
\|S(t) u\|_{L^p_t,L^q_x} \leq C \|u\|_{L^2}.
$$
\end{proposition}

Exploiting Bernstein's lemma and Littlewood-Paley theorem, we get the following.

\begin{corollary}\label{cor-StriBesov} Take $p,q \in [2,\infty[$ and $\sigma,\sigma_1 \geq 0$ such that
$$
\frac2{p} + \frac{d}{q} = \frac{d}{2} - \sigma_1,
$$
there exists $C,C'$ such that for all $u \in B_2^{\sigma_1, \sigma_1+\sigma}$,
$$
\|S(t)u\|_{L^p_t,W^{\sigma,q}}\leq C'\|S(t)u\|_{L^p_t,B_q^{0,\sigma}} \leq C \|u\|_{B_2^{\sigma_1, \sigma_1+\sigma}} \leq C \|u\|_{H^{\sigma_1+\sigma}}.
$$
\end{corollary}

\section{Principle} \lab{sec:principle}

In this section, we reduce the problem to finding a correct functional setting for our fixed point problem. \\

\noindent Writing $X = Y + Z$, the perturbation $Z$ satisfies
$$
i\partial_t Z = (m-\lap)Z + w*(2\Re \E(\overline Y Z) + |Z|^2)(Y+Z).
$$

\noindent Let an initial perturbation $Z_0 \in H^s$, we have that $Z$ solves the Cauchy problem
\begin{equation}\label{Cauchyprob}
\left \lbrace{\begin{array}{c} 
i\partial_t Z = (m-\lap)Z + w*(2\Re \E(\overline Y Z) + |Z|^2)(Y+Z)\\
Z_{|t=0} = Z_0
\end{array}} \right.
\end{equation}
if and only if the couple perturbation/induced potential $(Z,V)$ solves the Cauchy problem
$$
\left \lbrace{\begin{array}{c} 
i\partial_t Z = (m-\lap)Z + w*V (Y+Z),\\
V = 2\Re(\E (\bar Y Z)) + \E(|Z|^2) ,\\
Z_{|t=0} = Z_0.
\end{array}} \right.
$$

\noindent The idea is to set up spaces for $Z$ and $V$, $\Theta_Z$ and $\Theta_V$ such that $\Theta_Z \times \Theta_Z$ is embedded in $\Theta_V$ in the sense that there exists a constant $C$ such that for all $u,v \in \Theta_Z$,
$$
\|\E(uv) \|_{\Theta_V} \leq C \|u\|_{\Theta_Z}\|v\|_{\Theta_Z}.
$$
Indeed, one key idea behind the proof is that even if $V$ contains a linear term in $Z$, it behaves like a quadratic term on $Z$ in terms of functional spaces, which is better regarding the use of dispersion estimates, and this because of cancellations due to the randomness of the equation. \\

\noindent Using the Duhamel formulation of the equation on $Z$, we get that the Cauchy problem is equivalent to 
\be \lab{id:Zt}
Z(t) = S(t)Z_0 -i\int_{0}^t S(t-s) (w*V(s)) Z(s) ds -i\int_{0}^t S(t-s) (w*V(s)) Y(s) ds
\ee
and 
\begin{multline*}
V(t) = \E(|Z|^2)+2\Re \E(\bar Y(t) S(t)Z_0) - 2 \Re \E(i\bar Y(t) \int_{0}^t S(t-s) (w*V(s)) Z(s) ds \\
- 2 \Re \E(i\bar Y(t) \int_{0}^t S(t-s) (w*V(s)) Y(s) ds )
\end{multline*}
with the Schr\"odinger group $S(t) = e^{-it(m-\lap)}$.  \\

\noindent We set for all $V$, and all $Z$,
$$
W_V(Z) = -i\int_{0}^t S(t-s) (w*V(s)) Z(s) ds.
$$

\noindent In other terms, to solve the Cauchy problem, we have to solve the fix point
\be \lab{eq:pointfixe}
\mat{Z}{V} = A_{Z_0} \mat{Z}{V} = \mat{A_{Z_0}^1 \mat{Z}{V} }{ A^2_{Z_0}\mat{Z}{V}}
\ee
with 
$$
A_{Z_0}^1\mat{Z}{V} = S(t)Z_0 + W_V(Y) + W_V(Z)
$$
and
$$
A^2_{Z_0}\mat{Z}{V} = \E(|Z|^2) + 2\Re \E(\bar Y S(t) Z_0) + 2\Re \E( \bar Y W_V(Y)) + 2\Re \E(\bar Y(t) W_V(Z)).
$$

\noindent The map $A_{Z_0}$ has a constant part (or source term) $C_{Z_0}$ a linear part $L$ and a quadratic part $Q$ given by
\be \lab{def:az0}
A_{Z_0}\mat{Z}{V} = C_{Z_0} + L \mat{Z}{V} + Q\mat{Z}{V}
\ee
where
\be \lab{id:cz0}
C_{Z_0} = \mat{S(t) Z_0}{ 2\Re \E(\bar Y S(t) Z_0)},
\ee
\be \lab{id:L}
L\mat{Z}{V} = \mat{ W_V(Y)}{ 2\Re \E( \bar Y W_V(Y))},
\ee
and
\be \lab{id:Q}
Q\mat{Z}{V} = \mat{W_V(Z)}{ \E(|Z|^2)+ 2\Re \E(\bar Y W_V(Z))}.
\ee
The linear term can be written under the form
$$
L = \begin{pmatrix} 0 & L_2 \\ 0& L_1 \end{pmatrix}
$$
with $L_2(V) = W_V(Y)$ and $L_1(V) = 2\Re \E(\bar Y W_V(Y))$. \\

\noindent Solving the fixed point equation is now a problem of finding the right functional spaces and showing suitable continuity estimates. The properties we are going to show are summarised in the following proposition.

\begin{proposition}\label{prop-assumption} Let $\Theta = \Theta_Z \times \Theta_V$ and $\Theta_{0}$ be two Banach spaces such that
\begin{enumerate}
\item $\Theta_Z \times \Theta_Z$ is embedded in $\Theta_V$ in the sense that for all $u,v \in \Theta_Z$, $\|\E(uv)\|_{\Theta_V}\lesssim \|u\|_{\Theta_Z}\|v\|_{\Theta_Z}$ ,
\item there exists $C$ such that for all $Z_0 \in \Theta_{0} $, $\|C_{Z_0}\|_\Theta \leq C \|Z_0\|_{\Theta_{0}}$,
\item $1-L$ is continuous, invertible with continuous inverse as a linear operator of $\Theta$,
\item there exists $C_1$ such that for all $(Z,V ) \in \Theta$, $\|W_V(Z)\|_{\Theta_Z} \leq C_1 \|V\|_{\Theta_V}\|Z\|_{\Theta_Z}$,
\item there exists $C_2$ such that for all $(Z,V ) \in \Theta$, $\|2\Re \E(\bar Y W_V(Z))\|_{\Theta_V} \leq C_2 \|V\|_{\Theta_V}\|Z\|_{\Theta_Z}$,
\end{enumerate}
then there exists $\varepsilon > 0$ such that for all $Z_0 \in \Theta_{id}$ with $\|Z_0\|_{\Theta_{id}}\leq \varepsilon$, the Cauchy problem \eqref{Cauchyprob} has a unique solution in $\Theta_Z$ such that $2\Re \E(\bar Y Z) + \E(|Z|^2) \in \Theta_V$ and the flow thus defined is continuous in the initial datum.
\end{proposition}

The proof of this proposition and of the fact that finding such spaces implies the theorem is exposed in Section \ref{sec-proofth}. One important feature is that the terms $W_V(Z)$, $2\Re (\E(\bar Y W_V(Z)))$ and $\E (|Z|^2)$ are both bilinear and thus we can do a contraction argument directly from this setting.

\vspace{0.5cm}

\noindent We end this section by defining the solution spaces $\Theta_Z$ and $\Theta_V$.


\begin{definition}
Let for $d\geq 3$,
\begin{multline*}
\Theta_Z= \mathcal C(\R, H^s(\R^d,L^2_\omega)) \cap L^p(\R,W^{s,p}(\R^d,L^2_\omega)) \\
\cap L^{d+2}(\R,L^{d+2}(\R^d,L^2_\omega)) \cap L^4(\R,B_q^{0,\frac14}(\R^d,L^2_\omega))
\end{multline*}
where $s = \frac{d}{2} - 1$ is the critical regularity for the cubic Schr\"odinger equation, $p = 2 \frac{d+2}{d}$, $q= \frac{4d}{d+1}$.

Let
\be \lab{def:thetav}
\Theta_V = L^{\frac{d+2}{2}}_t, L^{\frac{d+2}{2}}_x + L^2_t,B_2^{-1/2,0}.
\ee

\end{definition}

In the next sections, we check that $\Theta = \Theta_Z \times \Theta_V$ satisfies assumptions 1,3,4, in Proposition \ref{prop-assumption} for $d\geq 3$ and assumption 5 only for $d\geq 4$.

\section{Embeddings and Strichartz estimates} \lab{sec:embeddings}

In this section, we check assumption 1 in Proposition \ref{prop-assumption}, and dispersion estimates for the linear flow which induce assumption 4. In the whole section, $d\geq 3$.

\begin{proposition} \lab{pr:contquadraZ} The space $ \Theta_Z \times \Theta_Z $ is embedded in $\Theta_V$ as in for all $u,v\in \Theta_Z$,
$$
\|\E(uv)\|_{\Theta_V} \lesssim \|u\|_{\Theta_Z} \|v\|_{\Theta_Z}.
$$
\end{proposition}

\begin{proof} We have that $L^{d+2}_t,L^{d+2}_x\times L^{d+2}_t,L^{d+2}_x$ is embedded in $L^{(d+2)/2}_t,L^{(d+2)/2}_x$.

It remains to prove that $L^4,B_q^{0,\frac14} \times L^4,B_q^{0,\frac14} $ is embedded in $L^2,B_2^{-1/2,0}$. The temporal part of the norm works by H\"older inequality. We are left with proving the embedding $B_q^{0,\frac14} \times B_q^{0,\frac14}$ in $B_2^{-1/2,0}$.

We use Lemma 4.1 in \cite{GNT3}, that we mentionned in the toolbox, Proposition \ref{prop-bilinBesov} with $s_1 = s_2 = 0$, $s_3 = \frac12$, $t_1=t_2= \frac14$, $t_3 = 0$, and $p_1 = p_2 = q$, $p_3=2$. We have for all $j = 1,2,3$, $p_j\geq 2$, $s_j < \frac{d}{p_j}$, $t_j \leq t_1+t_2+t_3$. We also have 
$$
\max(0,s_1,s_2,s_3,s_1+s_2+s_3) = \frac12 = d( \frac1{p_1} + \frac1{p_2}+ \frac1{p_3}-1)= t_1+t_2+t_3.
$$
Indeed, 
$$
d( \frac1{p_1} + \frac1{p_2}+ \frac1{p_3}-1) = d(\frac2{q} -\frac12) = d(\frac{d+1}{2d} - \frac12) = \frac12.
$$\end{proof}

\begin{proposition} There exists $C$ such that for all $u \in L^2_\omega,H^s$,
\be \lab{bd:cz01}
\|S(t) u\|_{\Theta_Z} \leq C \|u\|_{L^2_\omega,H^s}.
\ee
\end{proposition}

\begin{proof}First of all, let us mention that since $p,d+2$ and $q$ are bigger than $2$, we have, by Minkowski inequality for all $f \in L^2_\omega, \mathcal C(\R, H^s) \cap L^p,W^{s,p} \cap L^{d+2},L^{d+2}\cap L^4,B_q^{0,\frac14})$,
$$
\| f\|_{\Theta_Z} \leq \|f\|_{L^2_\omega, \mathcal C(\R, H^s) \cap L^2_\omega,L^p,W^{s,p} \cap L^2_\omega,L^{d+2},L^{d+2}\cap L^2_\omega,L^4,B_q^{0,\frac14}}.
$$

For $L^p,W^{s,p}$, we have
$$
\frac2{p}  + \frac{d}{p} = \frac{d+2}{p} = \frac{d}{2}
$$
hence Strichartz estimates, \ref{prop-Striadmiss}, apply.

For $L^{d+2},L^{d+2}$, we have that
$$
\frac2{d+2} + \frac{d}{d+2} = 1 = \frac{d}2 - s
$$
hence Strichartz inequalities, \ref{cor-StriBesov} ,apply.

For $L^4,B_q^{0,\frac14}$, we have
$$
\|S(t)u\|_{L^4_t ,B_q^{0,\frac14} } =\big\| \Big( \sum_{j<0} \|S(t)u_j\|_{L^q}^2 + \sum_{j\geq 0} 2^{j/2} \|S(t)u_j\|_{L^q}^2\Big)^{1/2}\big\|_{L^4_t}.
$$
Since $4\geq 2$, by convexity we have,
$$
\|S(t)u\|_{L^4_t ,B_q^{0,\frac14} } \leq \Big( \sum_{j<0} \|S(t)u_j\|_{L^4,L^q}^2 + \sum_{j\geq 0} 2^{j/2} \|S(t)u_j\|_{L^4,L^q}^2\Big)^{1/2}.
$$
Let $q_1 = \frac{2d}{d-1}$ and $s_1 = \frac{d-3}{4}$. Since $d\geq 3$ and $s_1 = \frac{d}{q_1}-\frac{d}{q}$ (we recall $q =\frac{4d}{d+1}$) , we have by Bernstein lemma,
$$
\|S(t)u_j\|_{L^q}\leq C 2^{js_1} \|S(t)u_j\|_{L^{q_1}}
$$
and since
$$
\frac12 + \frac{d}{q_1} = \frac{d}{2}
$$
we get by Strichartz estimates, \ref{prop-Striadmiss},
$$
\|S(t)u_j\|_{L^4_t,L^{q_1}_x} \leq C\|u_j\|_{L^2}.
$$
We deduce
$$
\|S(t)u\|_{L^4_tB_q^{0,\frac14} } \leq  C \|u\|_{B_2^{s_1,s_1+1/4}}.
$$
We have $s_1\geq 0$ and $s_1 + \frac14 = \frac{d-1}{4} \leq \frac{d-2}{2} = s$, by \ref{prop-Besovembed}, we get
$$
\|S(t)u\|_{L^4_t, B_q^{0,\frac14} } \leq  C \|u\|_{B_2^{0,s}}\leq C \|u\|_{H^s}.
$$

\end{proof}

\begin{proposition} \lab{pr:contWVZ} Let $d\geq 3$. There exists $C_1$ such that for all $Z\in \Theta_Z$ and $V\in \Theta_V$,
$$
\|W_V(Z)\|_{\Theta_Z} \leq C_1 \|Z\|_{\Theta_Z}\|V\|_{L^{(d+2)/2}_{t,x}}\leq C_1 \|Z\|_{\Theta_Z}\|V\|_{\Theta_V}.
$$
\end{proposition}

\begin{proof}

Thanks to the previous proposition, Christ-Kiselev lemma and dual Strichartz estimates, we have
$$
\|W_V(Z)\|_{\Theta_Z} \lesssim \|(w*V) u\|_{L^{p'},W^{s,p'},L^2_\omega}
$$
with $ p '$ the conjugate of $p$ that is $p'= 2\frac{d+2}{d+4}$. We have $\frac1{p'} = \frac1{p} + \frac2{d+2}$ hence, thanks to H\"older's inequality,
$$
\|W_V(u)\|_{\Theta_Z} \lesssim \|w*V\|_{W^{s,(d+2)/2}}\|Z\|_{L^{p},W^{s,p},L^2_\omega}.
$$

The potential $w$ allows us to lose $s$ derivatives, that is,
$$
\|w*V\|_{W^{s,(d+2)/2}} \leq \|w\|_{W^{s,1}} \|V\|_{L^{(d+2)/2}_{t,x}},
$$
and we get
$$
\|W_V(u)\|_{\Theta_Z} \lesssim \|V\|_{\Theta_V}\|Z\|_{\Theta_Z}.
$$

\end{proof}

\section{Linear term} \lab{sec:lin}

We study in this section the linearised operator $L$ defined in \fref{id:L}. We prove that $1-L$ is continuous, invertible with continuous inverse on $\Theta_Z\times \Theta_V$. We recall that
\be \lab{id:1-l}
1-L = \begin{pmatrix} 1 & -L_2 \\ 0 & 1-L_1 \end{pmatrix},
\ee
hence it is sufficient to prove that $1-L_1$ is continuous, invertible with continuous inverse from $\Theta_V$ to $\Theta_V$ and that $L_2 $ is continuous from $\Theta_V$ to $\Theta_Z$. We start with the continuity of $L_2$.

\begin{proposition}\label{prop-cont} The operator
$$
L_2 : V \mapsto W_V(Y)=-i\int_{0}^t S(t-s)\left[ (w*V(s)) Y(s) \right]ds.
$$
is continuous from $\Theta_V$ to $\Theta_Z$.
\end{proposition}

Proposition \ref{prop-cont} is a corollary of the following estimates.

\begin{proposition}\label{prop-contU} Let $\sigma,\sigma_1,p_1,q_1$ be such that $\sigma,\sigma_1 \geq 0$, $\sigma \leq \lceil s \rceil$, $p_1>2$ and
$$
\frac{2}{p_1} + \frac{d}{q_1}= \frac{d}{2}-\sigma_1.
$$
Assume moreover, that for all $\alpha \in \mathbb N^d$ with $|\alpha|\leq 2\lceil \sigma \rceil$, $|\pa^\alpha h|\lesssim \langle \xi \rangle^{-2}$. There exists $C>0$ such that for all $U \in L^2_t ,B_2^{-\frac12,\sigma_1+ \sigma}$, 
$$
\sup_{A\in \R} \big\|\int_{0}^A S(t-s) Y(s) U(s) ds\big\|_{L^{p_1}_t,B_{q_1}^{0,\sigma},L^2_\omega}  \leq C \|U\|_{L^2_t , B_2^{-\frac12,\sigma_1+\sigma}} 
$$
and if $\sigma + \sigma_1 \leq \lceil s \rceil$, there exists $C>0$ such that for all $U \in L^2_t, B_2^{\sigma_1-\frac12,\sigma_1+ \sigma} $,
$$
\sup_{A\in \R } \big\|\int_{0}^A S(t-s) Y(s) U(s) ds\big\|_{L^{p_1}_t ,W^{\sigma,q_1},L^2_\omega}\leq C \|U\|_{L^2_t ,B_2^{\sigma_1-\frac12,\sigma_1+\sigma}} .
$$
\end{proposition}

\begin{proof}[Proof of Proposition \ref{prop-cont}]Take $p_1>2$, $\sigma_1 \geq 0$ and $q_1\geq 2$ such that
$$
\frac2{p_1} + \frac{d}{q_1} = \frac{d}{2} - \sigma_1.
$$
Assuming Proposition \ref{prop-contU}, we get by Christ-Kiselev lemma, since $p_1>2$,
$$
\big\|\int_{0}^t S(t-s) Y(s) U(s) ds\big\|_{L^{p_1}_t, B_{q_1}^{0,\sigma} ,L^2_\omega} \leq C \|U\|_{L^2_t, B_2^{-\frac12,\sigma_1+\sigma} }.
$$
We use it with $p_1 = 4$, $q_1 = q$, $\sigma = \frac14$ and $\sigma_1 = \frac{d-3}{4} \leq s-\sigma$, to get
$$
\big\|\int_{0}^t S(t-s) Y(s) U(s) ds\big\|_{L^{4}_t, B_{q}^{0,1/4} ,L^2_\omega} \leq C \|U\|_{L^2_t ,B_2^{-\frac12,s} }.
$$
We also have 
$$
\big\|\int_{0}^t S(t-s) Y(s) U(s) ds\big\|_{L^{p_1}_t ,W^{q_1,\sigma} ,L^2_\omega} \leq C \|U\|_{L^2_t ,B_2^{\sigma_1-\frac12,\sigma_1+\sigma}}.
$$
We apply it to $(p_1,q_1,\sigma_1,\sigma)$ equal to either
$$
(p,p,0,s), \; (d+2,d+2,s,0) \textrm{ or } (\infty,2,0,s)
$$
to get
$$
\|W_V(Y)\|_{\Theta_Z} \leq C \|w*V\|_{L^2_t, B_2^{-1/2,s} }
$$
and thus, by putting the derivatives on high frequency on $w$ and on low frequency on $V$,
$$
\|W_V(Y)\|_{\Theta_Z} \leq C \|w\|_{W^{s,1}}\|V\|_{L^2_t, B_2^{-1/2,0}}.
$$
\end{proof}

It remains to prove Proposition \ref{prop-contU}. To do so, we first establish preliminary lemmas.

\begin{lemma} Let $U\in L^{\infty}_t L^2_x$ and define $L_3^A(U) = \int_{0}^A S(t-s) Y(s) U(s) ds$. We have
\be \lab{id:L3A}
\E(|L_3^A(U)|^2) = \int d\eta |f(\eta)|^2 \Big| \int_{0}^A S_\eta(t-s) U(s) ds \Big|^2,
\ee
where
\be \lab{eq:defSeta}
S_\eta(t) = e^{-it(-\lap - 2i\eta\cdot \grad)}.
\ee
\end{lemma}

\begin{proof}

Let $U\in L^{\infty}_t L^2_x$. We write $U(s)Y(s)$ as a Wiener integral:
$$
U(s)Y(s)=\int_{\eta \in \mathbb R^d} f(\eta)e^{i\eta.x-is(m+|\eta|^2)} U(s) dW(\eta).
$$
Above, we remark that thanks to Fubini and to the Wiener integration
\bee
\|U(s)Y(s)\|_{L^2_\omega, L^2_x}^2 & = & \int \E \left(\left|\int_{\eta \in \mathbb R^d} f(\eta)e^{i\eta.x-is(m+|\eta|^2)} U(s,x) dW(\eta)\right|\right)dx \\
& = & \int |U(s,x)|^2\int_{\eta \in \mathbb R^d} |f(\eta)|^2d\eta dx= \| U(s)\|_{L^2_x}^2 \| f \|_{L^2}^2.
\eee
Therefore, almost everywhere $U(s)Y(s,\omega)\in L^2_x$, making $S(t-s)(U(s)Y(s))$ well defined. From the commutator relation $S(t)(e^{i\eta.x}U)=e^{i\eta.x- it(m+|\eta|^2)}S_\eta(t)U$ we infer that:
\bea
\non S(t-s)(U(s)Y(s)) & = & \int_{\eta \in \mathbb R^d} f(\eta)e^{i\eta.x-is(m+|\eta|^2)-i(t-s)(m+|\eta|^2)} S_\eta (t-s)(U(s)) dW(\eta) \\
\lab{id:YU} & = & \int_{\eta \in \mathbb R^d} f(\eta)e^{i\eta.x-it (m+ |\eta|^2)} S_\eta (t-s)(U(s)) dW(\eta).
\eea
From the definition of the Wiener integral, taking the expectation one obtains:
\bee
 \E \left( \left| \int_0^A S(t-s)(U(s)Y(s))ds\right|^2 \right) & = & \E \Bigl( \int_{\eta \in \mathbb R^d}\int_0^A  f(\eta)e^{i\eta.x-it (m+ |\eta|^2)} S_\eta (t-s)(U(s)) dW(\eta)ds\\
&& \times \int_{\eta' \in \mathbb R^d} \int_0^A  \bar f(\eta ')e^{-i\eta'.x+it (m+ |\eta'|^2)} \overline{S_{\eta'} (t-s')(U(s'))} dW(\eta ')ds' \Bigr)\\
&=& \int_{\eta \in \mathbb R^d} \int_0^A \int_0^A |f(\eta)|^2 S_{\eta}(t-s)(U(s)) \overline{S_{\eta} (t-s')(U(s'))} ds ds' d\eta \\
&=& \int_{\eta \in \mathbb R^d} |f(\eta)|^2 \Big| \int_{0}^A S_\eta(t-s) U(s) ds \Big|^2d\eta.
\eee
which ends the proof of the identity \fref{id:L3A}.

\end{proof}

We start the proof of Proposition \ref{prop-contU} with the case $\sigma=0$ for H\"older spaces and $\sigma=\sigma_1 = 0$ for Besov spaces.

\begin{lemma} \lab{lem:l3a}

Let $\sigma_1,p_1,q_1$ be such that $\sigma_1 \geq 0$, $p_1>2$ and
$$
\frac{2}{p_1} + \frac{d}{q_1}= \frac{d}{2}-\sigma_1.
$$
There exists $C>0$ such that for all $U \in L^2_t, B_2^{\sigma_1-\frac12,\sigma_1+\sigma}$, 
\be \lab{bd:leml3a}
\sup_{A\in \R} \|L_3^A(U)\|_{L^{p_1}_t, L^{q_1}_x ,L^2_\omega} \leq C \|U\|_{L^2_t ,B_2^{\sigma_1-\frac12,\sigma_1}}.
\ee
\end{lemma}

\begin{proof} We start by taking $U$ in the Schwartz class to allow the following computations and we conclude by density. We have from \fref{id:L3A} and Minkowski's inequality:
\bee
\|L_3^A(U)\|_{L^{p_1}_t ,L^{q_1}_x ,L^2_\omega}^2 &=& \big\| \int_{\mathbb R^d} d\eta |f(\eta)|^2 \Big| \int_{0}^A S_\eta(t-s) U(s) ds \Big|^2 \big\|_{L^{p_1/2}_t ,L^{q_1/2}_x} \\
&\leq & \int_{ \mathbb R^d} d\eta |f(\eta)|^2  \big\|  \int_{0}^A S_\eta(t-s) U(s) ds  \big\|_{L^{p_1}_t,L^{q_1}_x}^2.
\eee
By Strichartz inequality and Bernstein lemma, Corollary \ref{cor-StriBesov}, we get
$$
\|L_3^A(U)\|_{L^{p_1}_t ,L^{q_1}_x ,L^2_\omega}^2 \leq \int_{\mathbb R^d}  d\eta |f(\eta)|^2 \big\|  \int_{0}^A S_\eta(-s) U(s) ds  \big\|_{B_2^{\sigma_1,\sigma_1}}^2.
$$
We introduce $U_1$ defined by $\hat U_1(\xi) = |\xi|^{\sigma_1} \hat U(\xi)$, we have 
$$
\big\|  \int_{0}^A S_\eta(-s) U(s) ds  \big\|_{B_2^{\sigma_1,\sigma_1}}^2 = \int_{\mathbb R^d} d\xi \int_{0}^A dt_1 \int_{0}^A dt_2 e^{i(t_1-t_2)(\xi^2 - 2\xi\cdot \eta)} \overline{\hat U_1 (\xi,t_1)} \hat U_1(\xi,t_2).
$$
We do the change of variables $t= t_2-t_1$, we get
$$
\big\|  \int_{0}^A S_\eta(-s) U(s) ds  \big\|_{B_2^{\sigma_1,\sigma_1}}^2 = \int_{\mathbb R^d} d\xi \int_{\R} dt \int_{t_1\in D_t} dt_1 e^{-it(\xi^2 - 2\xi\cdot \eta)} \overline{\hat U_1 (\xi,t_1)} \hat U_1(\xi,t+t_1),
$$
where $D_t=[-t,A-t]\cup [0,A]$. We integrate over $\eta$ reminding that $h$ is the inverse Fourier transform of $|f|^2$, we get
$$
\|L_3^A(U)\|_{L^{p_1}_t, L^{q_1}_x,L^2_\omega}^2\lesssim \int_{\mathbb R^d} d\xi \int_{\R} dt \int_{t_1\in D_t} dt_1 e^{-it\xi^2 } h(-2t\xi)\overline{\hat U_1 (\xi,t_1)} \hat U_1(\xi,t+t_1).
$$
We use Cauchy-Schwarz inequality over $t_1$ to get
$$
\|L_3^A(U)\|_{L^{p_1}_t, L^{q_1}_x, L^2_\omega}^2\lesssim \int_{\mathbb R^d} d\xi \int_{\R} dt | h(-2t\xi)|\|\hat U_1 (\xi,t_1)\|_{L^2_{t_1}}^2 
$$
We have $|h(-2t\xi)| \leq \an{t|\xi|}^{-2}$ from the hypotheses of Theorem \ref{thmmain}. We then do the change of variable $\tau = t|\xi|$ and get
$$
\|L_3^A(U)\|_{L^{p_1}_t, L^{q_1}_x ,L^2_\omega}^2\lesssim \int_{\mathbb R^d} d\xi \int_{\R} d\tau \an{\tau}^{-2} |\xi|^{-1}\|\hat U_1 (\xi,t_1)\|_{L^2_{t_1}}^2
$$
and since $\an{\tau}^{-2}$ is integrable, we have
$$
\|L_3^A(U)\|_{L^{p_1}_t,L^{q_1}_x,L^2_\omega}^2\lesssim \int_{ \mathbb R^d} \| |\xi|^{-1/2}\hat U_1 (\xi,\cdot)\|_{L^2_t}^2d\xi = \|U\|_{L^2_t ,B_2^{\sigma_1-1/2,\sigma_1-1/2}}^2
$$
which is the desired result.

\end{proof}

\begin{lemma} Let $p_1,q_1$ be such that $p_1>2$ and
$$
\frac{2}{p_1} + \frac{d}{q_1}= \frac{d}{2}.
$$
There exists $C>0$ such that for all $U \in L^2_t ,B_2^{-\frac12,0}$, 
\be \lab{bd:leml3a2}
\sup_{A\in \R} \|L_3^A(U)\|_{L^{p_1}_t ,B_{q_1}^{0,0},L^2_\omega} \leq C \|U\|_{L^2_t ,B_2^{-\frac12,0}}.
\ee
\end{lemma}

\begin{proof} By Minskowski's inequality, we have 
$$
\|L_3^A(U)\|_{L^{p_1}_t ,B_{q_1}^{0,0},L^2_\omega} \leq \|L_3^A(U)\|_{L^2_\omega L^{p_1}_t,B_{q_1}^{0,0}} .
$$
By Strichartz inequality we get
$$
\|L_3^A(U)\|_{L^{p_1}_t, B_{q_1}^{0,0},L^2_\omega} \lesssim \big\|\int_{0}^AS(-s) Y(s) U(s) ds\|_{L^2_\omega ,L^2_x} .
$$
We use the formula \fref{id:L3A} for $L_3^A(U)$ when $t=0$ and obtain
$$
\|L_3^A(U)\|_{L^{p_1},B_{q_1}^{0,0},L^2_\omega}^2 \lesssim \int |f(\eta)|^2 \big\|\int_{0}^AS_\eta (-s)  U(s) ds\|_{L^2_\omega,L^2}.
$$
Following exactly the same steps as in the proof of Lemma \ref{lem:l3a}, we obtain:
$$
\|L_3^A(U)\|_{L^{p_1},B_{q_1}^{0,0},L^2_\omega}^2 \lesssim \| U \|_{B_2^{-\frac 12,-\frac 12}}\lesssim \| U \|_{B_2^{-\frac 12,0}},
$$
which ends the proof of the lemma.

\end{proof}

We can now end the proof of Proposition \ref{prop-contU} thanks to the three lemmas above.

\begin{proof}[Proof of proposition \ref{prop-contU}.]

 For vectors $\alpha  = (\alpha_1,\hdots ,\alpha_d)\in \N^d$ and $\eta = (\eta_1,\hdots ,\eta_d) \in \R^d$ we write $|\alpha | = \sum_{j=1}^d \alpha_j$,  $\partial^\alpha = \prod_{j=1}^d \partial_j^{\alpha_j}$ and
$$
C_\alpha^\beta  = \prod_{j=1}^d C_{\alpha_j}^{\beta_j} \textrm{ and } \eta^\alpha = \prod_j \eta_j^{\alpha_j} 
$$
where $C_{\alpha_j}^{\beta_j}$ is a binomial coefficient. We have for all $\alpha \in \N^d$,
$$
\partial^\alpha L_3^A(U) = \sum_{\beta + \gamma = \alpha}C_\alpha^\beta \int_{0}^A S(t-s)  \partial^\beta Y(s) \partial^\gamma U(s)ds.
$$
Indeed, $Y$ is almost everywhere differentiable and there holds, for $|\beta| \leq \lceil s\rceil$ :
$$
\partial^\beta Y(s)= \int_{ \mathbb R^d} i^{|\beta|}\eta^\beta f(\eta) e^{i \eta.x-is(m+|\eta|^2)}dW(\eta),
$$
meaning that replacing $Y$ by $\partial^\beta Y$ consists in replacing $f(\eta)$ by $i^{|\beta|}\eta^\beta f(\eta)$. Therefore, we have for $\sigma_1 \geq 0$, $p_1,q_1>2$ such that $\frac2{p_1} + \frac{d}{q_1} = \frac{d}{2} - \sigma_1$, and $\sigma \in \N \cap [0, \lceil s \rceil]$, for any $A\in \mathbb R$:
$$
\|L_3^A(U)\|_{L^{p_1}_t,W^{\sigma,q_1},L^2_\omega} \leq C \|U\|_{L^2_t,B_2^{\sigma_1-1/2, \sigma_1+\sigma}},
$$
thanks to Lemma \ref{bd:leml3a}. Above, the constant $C$ depends on $\sup_{\alpha \in \mathbb N^d, \ |\alpha|\leq 2 \sigma } \| \langle \xi \rangle^2 \pa^{\alpha} h \|_{L^{\infty}}$, because the Fourier transform of $|\eta^{\alpha}|^2|f|^2$ is $\pa^{2\alpha}h$ up to multiplication by a constant. By interpolation,
$$
\|L_3^A(U)\|_{L^{p_1}_t,W^{\sigma,q_1},L^2_\omega} \leq C \|U\|_{L^2_t,B_2^{\sigma_1-1/2, \sigma_1+\sigma}}
$$
for all $\sigma \geq 0$. And for $\sigma + \sigma_1\in \N$ and $\frac2{p_1} + \frac{d}{q_1} = \frac{d}{2} - \sigma_1$, we have by Bernstein lemma
$$
\|L_3^A(U)\|_{L^{p_1},B_{q_1}^{0,\sigma},L^2_\omega} \lesssim \|L_3^A(U)\|_{L^{p_1},B_{q_2}^{0,\sigma+\sigma_1},L^2_\omega}
$$
with $\frac2{p_1} + \frac{d}{q_2} = \frac{d}{2} $. Therefore, applying Lemma \ref{bd:leml3a2} with $\sigma + \sigma_1 \leq \lceil s \rceil$:
$$
\|L_3^A(U)\|_{L^{p_1}_t,B_{q_1}^{0,\sigma},L^2_\omega} \leq C \|U\|_{L^2_t,B_2^{-1/2, \sigma_1+\sigma}}.
$$
We get the result for $\sigma \geq 0, \sigma_1 \geq 0, \sigma + \sigma_1 \leq \lceil s \rceil$ by interpolation.
\end{proof}

Having established the desired estimates for $L_2$, we now deal with the other part of the linear term, which involves $L_1$. We start by computing an explicit formula for this operator, and then show some continuity properties.

\begin{lemma}\label{lem-l1} The operator $L_1$ is a Fourier multiplier of symbol $\hat w m_f$ (in both space and time) given by
for all $t\in \R$ and all $\xi \in \R^d$,
$$
\mathcal F_x\left(L_1(V)\right)(t,\xi) = -2 \hat w(\xi)\int_{0}^t \sin(|\xi|^2(t-s)) h(2\xi(t-s)) \mathcal F_x V(s,\xi) ds,
$$
where $h$ is the inverse Fourier transform of $|f|^2$, or, put another way:
$$
\mathcal F_{t,x}\left(L_1(V)\right)(\tau,\xi) = \hat w (\xi)m_f(\tau,\xi) \mathcal F_{t,x} V(\tau,\xi),
$$
where 
\be \lab{def:mf}
m_f(\tau,\xi)= -2\mathcal F_t \left(\sin(|\xi|^2t) h(2\xi t)1_{t\geq 0} \right)(\tau)=-2 \int_0^{+\infty} e^{-i\tau t} \sin (|\xi|^2t)h(2\xi t) dt .
\ee

\end{lemma}

\begin{proof}[ Proof of Lemma \ref{lem-l1}.] Let $V\in L^{\infty}_tL^2_x$ so that one can perform the computations below. Using the formula \fref{id:YU} we write $W_V (Y)$ as a Wiener integral:
$$
W_V (Y)=-i\int_{\mathbb R^d} \int_0^t f(\eta) e^{i\eta.x-it(m+|\eta|^2)}S_\eta (t-s)(w*V(s))dsdW(\eta).
$$
By the property of Wiener integration:
\bee
\E \left(W_V (Y)\bar Y\right)&=&-i\E \Bigl( \int_{\mathbb R^d} \int_0^t f(\eta) e^{i\eta.x-it(m+|\eta|^2)}S_\eta (t-s)(w*V(s))dsdW(\eta)\\
&& \times \int_{\mathbb R^d} \bar f (\eta')e^{-i\eta'.x+it(m+|\eta'|^2)}d\bar W(\eta')\Bigr) \\
&=&-i  \int_{\mathbb R^d} \int_0^t |f(\eta)|^2 S_\eta (t-s)(w*V(s))dsd \eta.
\eee
Set $T(t,\xi)$ the Fourier transform in space of $\E\Big( W_V(Y)\overline{ Y }\Big)$. From \fref{eq:defSeta} and Fubini we infer:
\bee
T(t,\xi)&=& -i  \int_{\mathbb R^d} \int_0^t |f(\eta)|^2 e^{-i(t-s)(|\xi|^2+2\xi.\eta)} \hat w(\xi) \hat V(s,\xi)dsd \eta\\
&=& -i \int_0^t h(2\xi (t-s)) e^{-i(t-s)|\xi|^2} \hat w(\xi) \hat V(s,\xi)dsd \eta.
\eee
The space Fourier transform of $2\Re \E (W_V^1(Y)\overline Y)$ is $T(t,\xi)+\overline{T(t,-\xi)}$. Given that $w*V$ is real and $|f|^2$ is even and real, and thus $h$ is real and even, we have
$$
T(t,\xi) + \overline T(t,-\xi) = -2 \int_{0}^t \sin( \xi^2 (t-s)) h(2\xi(t-s))\hat w(\xi) \hat V(s,\xi) ds
$$
which gives the result. 

\end{proof}

Since $ L_1$ is a Fourier multiplier with symbol $\hat w m_f$ given by \fref{def:mf}, the continuity and the invertibility of $1- L_1$ on $L^2$ based space-time functional spaces is equivalent to the boundedness and the non-vanishing of $1-\hat w m_f$. This issue was studied in \cite{LS2} where the author showed the following:

\begin{proposition}[Lewin Sabin,\cite{LS2}] \lab{pr:bd 1-L1 L2}

Let $d\geq 1$, Let $d\geq 4$. Let $s = \frac{d}2 - 1$. Let $f$ be a radial map in $L^2\cap L^\infty(\R^d)$. Consider the Fourier multiplier $m_f$ defined by \fref{def:mf}.
Assume : \begin{itemize} 
\item $\int_{\R^d} |\xi|^{1-d}|f\grad  f| < \infty$,
\item writing $r= |\xi|$ the radial variable, $\partial_r |f|^2 <0$ for $r>0$,
\item $|\xi|^{2-d}h \in L^1(\R^d)$.
\end{itemize}
If $w\in L^1(\mathbb R^d)$ is an even function such that
$$
\| (\hat w)_-\|_{L^{\infty}}\left(\int_{\mathbb R^d} \frac{|h|}{|x|^{d-2}}dx\right)<2|\mathbb S^{d-1}|,
$$
and such that
$$
\epsilon_g \hat w (0)_+<2|\mathbb S^{d-1}|, \ \ \text{where} \ \ \epsilon_g:=-\liminf_{(\tau,\xi)\rightarrow (0,0)} \frac{\Re m_f(\tau,\xi)}{2|\mathbb S^{d-1}|},
$$
then there holds
\be \lab{bd:minhatwmf}
\min_{(\omega,\xi)\in \mathbb R \times \mathbb R^d} |\hat w(k)||m_f(\tau,k)-1|>0
\ee
and the operator $1-L_1$ is invertible on $L^2_{t,x}$ and $L^2_tB^{-\frac 12,0}_2$.

\end{proposition}

\begin{remark}
The $f$ in the paper by Lewin and Sabin corresponds to $g(r)= |f(\sqrt r \vec e)|^2$ for any unitary vector of $\R^d$, $\vec e$. The assumptions they make on $g$ are implied by the ones we make on $f$.\end{remark}

\begin{proof}

For the proof, we refer to \cite{LS2}. There the authors show that under the assumptions of the proposition, $(\tau,\xi)\mapsto \hat w(\xi) m_f(\tau,\xi)$ is uniformly bounded on $\mathbb R \times \mathbb R^d$, and that \fref{bd:minhatwmf} holds. Therefore, $(1-\hat w(k)m_f(\tau,k))^{-1}$ is a bounded function. This implies the continuity of $1-L_1$ both on $L^2_{t,x}$ and $L^2_tB^{-\frac 12,0}_2$.

\end{proof}

The continuity of $1-L_1$ on $\Theta_V$ is a consequence of the above continuity on $L^2_{t,x}$ and of other mild assumptions on $f$ and $w$.

\begin{proposition}[Lewin-Sabin, \cite{LS2}] \lab{pr:cont1-L1}

Assume that the hypothesis of Proposition \ref{pr:bd 1-L1 L2} hold and that moreover $\int (|h|+ |\nabla h|)dr<+\infty$ (that is $h$ and $\grad h$ in $|\xi|^{d-1} L^1(\R^d)$) and $(1+|\xi|)\hat w \in L^{2(d+2)/(d-2)}$, then $1-L_1$ and $(1-L_1)^{-1}$ are continuous from $\Theta_V$ into $L^{\frac{d+2}{2}}_{t,x}$.

\end{proposition}

\begin{proof}

A first direction to prove Proposition \ref{pr:cont1-L1} is to show that the Fourier multiplier $w(\xi)m_f(\tau,\xi)$ and its inverse are continuous operators on Lebesgue spaces using standard harmonic analysis tools. This is done in \cite{LS2} where the authors show that this multiplier satisfies suitable conditions ensuring the application of Stein's and Marcinkiewicz's theorems. One can check using the computations there that its inverse also satisfies the same suitable conditions. We give here another proof avoiding the use of advanced harmonic analysis. We claim that
\be \lab{bd:multiplier}
\hat w (\xi)m_f(\tau,\xi)\in L^{2\frac{d+2}{d-2}}(\mathbb R^{1+d}).
\ee
Assuming the above bound, then one computes the following, using the continuity of the Fourier transform from $L^p$ into $L^{p'}$ for $1\leq p<+\infty$ and H\"older inequality. We have by definition
$$
\| L_1 V \|_{L^{\frac{d+2}{2}}_{t,x}}  \lesssim  \| \hat w (\xi)m_f(\tau,\xi) \mathcal F_{t,x}V \|_{L^{\frac{d+2}{d}}(\mathbb R^{1+d})}.
$$ 
We use H\"older inequality to get
$$
\| L_1 V \|_{L^{\frac{d+2}{2}}_{t,x}}  \lesssim  \| \hat w (\xi)m_f(\tau,\xi) \|_{L^{2\frac{d+2}{d-2}}(\mathbb R^{1+d})} \| \mathcal F_{t,x}V \|_{L^2(\mathbb R^{1+d})} 
$$
Because of the bound \eqref{bd:multiplier}, and that $\| \mathcal F_{t,x}V \|_{L^2(\mathbb R^{1+d})} \leq \| V \|_{L^2_t,B^{-\frac 12,0}}$, we have
$$
\| L_1 V \|_{L^{\frac{d+2}{2}}_{t,x}}  \lesssim  \| V \|_{L^2_t,B^{-\frac 12,0}} 
$$
which is in turn less than
$ \| V \|_{\Theta_V}$. This shows the continuity of $L_1$ from $\Theta_V$ into $L^{\frac{d+2}{2}}_{t,x}$. Similarly, using the fact that $(1-\hat w (\xi)m_f(\tau,\xi))^{-1}$ is uniformly bounded from Proposition \ref{pr:bd 1-L1 L2}, we get through algebraic considerations
$$
\| (1-L_1)^{-1} V \|_{L^{\frac{d+2}{2}}_{t,x}} \lesssim  \| V \|_{L^{\frac{d+2}{2}}_{t,x}} +\| (1-(1-L_1)^{-1})V \|_{L^{\frac{d+2}{2}}_{t,x}} .
$$
By definition, and using that for $p\geq 2$, and $\chi \in L^{p'}$, $\|\hat \chi\|_{L^p} \leq \|\chi\|_{L^{p'}}$, we have
$$
\| (1-(1-L_1)^{-1})V \|_{L^{\frac{d+2}{2}}_{t,x}} \leq  \| \frac{\hat w (\xi)m_f(\tau,\xi)}{1-\hat w (\xi)m_f(\tau,\xi)} \mathcal F_{t,x}V \|_{L^{\frac{d+2}{d}}(\mathbb R^{1+d})} .
$$
Using that $1-\hat w (\xi)m_f(\tau,\xi)$ is bounded from bellow we get
$$
\| (1-(1-L_1)^{-1})V \|_{L^{\frac{d+2}{2}}_{t,x}} \lesssim \| \hat w (\xi)m_f(\tau,\xi) \|_{L^{2\frac{d+2}{d-2}}(\mathbb R^{1+d})} \| \mathcal F_{t,x}V \|_{L^2(\mathbb R^{1+d})} 
$$
and we conclude as previously. This shows the continuity of $(1-L_1)^{-1}$ from $\Theta_V$ into $L^{\frac{d+2}{2}}_{t,x}$. Hence it remains to prove \eqref{bd:multiplier}. We compute for $|\xi|>1$ and $\tau \neq 0$, doing the change of variables $t\leftarrow t|\xi|$ :
$$
m_f(\tau,\xi)  =  -2\int_0^{+\infty} e^{-i\tau t} \sin (|\xi|^2t)h(2\xi t)dt=-\frac{2}{|\xi|} \int_0^{+\infty} e^{-i\frac{\tau}{|\xi|} t} \sin (|\xi|t)h\left(2\frac{\xi}{|\xi|} t\right)dt .
$$
We write $\frac1{|\xi|} e^{-i\frac{\tau}{|\xi|}t} = i\tau^{-1}\partial_t e^{-i\frac{\tau}{|\xi|}t}$ to get
$$
m_f(\tau,\xi) = -\frac{2i}{\tau} \int_0^{+\infty} \pa_t \left(e^{-i\frac{\tau}{|\xi|} t}\right) \sin (|\xi|t)h\left(2\frac{\xi}{|\xi|} t\right)dt 
$$
and we integrate by parts to get
$$
m_f(\tau,\xi)
=\frac{2i}{\tau} \left( |\xi| \int_0^{+\infty} e^{-i\frac{\tau}{|\xi|} t} \cos (|\xi|t)h\left(2\frac{\xi}{|\xi|} t\right)dt+2\int_0^{+\infty} e^{-i\frac{\tau}{|\xi|} t} \sin (|\xi|t)\frac{\xi}{|\xi|}.\nabla h\left(2\frac{\xi}{|\xi|} t\right)dt   \right) 
$$
Finally, we use the fact that $h$ is radially symmetric and the integrability assumptions to get
$$
m_f(\tau,\xi)
\lesssim  \frac{1+|\xi|}{|\tau|}
$$
Since $m_f$ is uniformly bounded from Proposition \ref{pr:bd 1-L1 L2}, one deduces that $|m_f(\tau,\xi) |\lesssim (1+|\xi|)(1+|\tau|)^{-1}$. Therefore one concludes that:
$$
\int_{\mathbb R^{1+d}} |\hat w (\xi)m_f(\tau,\xi)|^{2\frac{d+2}{d-2}} d\tau d\xi\lesssim \int_{\mathbb R^{1+d}} (1+|\xi|)^{2\frac{d+2}{d-2}}|\hat w|^{2\frac{d+2}{d-2}}(\xi)\frac{1}{(1+|\tau|)^{2\frac{d+2}{d-2}}}d\xi d \tau
$$
and since $\frac{1}{(1+|\tau|)^{2\frac{d+2}{d-2}}}$ is integrable, we have
$$
\int_{\mathbb R^{1+d}} |\hat w (\xi)m_f(\tau,\xi)|^{2\frac{d+2}{d-2}}\lesssim \| (1+|\xi|)\hat w\|_{L^{2\frac{d+2}{d-2}}(\mathbb R^d)}<+\infty.
$$

\end{proof}


\section{The remaining quadratic term} \lab{sec:quad}

We have already dealt with the quadratic terms $\E(|Z|^2)$ and $W_V(Z)$ in Propositions \ref{pr:contquadraZ} and \ref{pr:contWVZ} respectively, it remains to prove assumption 5 in Proposition \ref{prop-assumption}. In what follows, we assume $d\geq 4$, it is the only part of the paper that does not work in dimension $3$.

\begin{proposition} \lab{pr:contEWVZ} There exists $C>0$ such that for all $(Z,V ) \in \Theta$, 
$$
\|2\Re \E(\bar Y W_V(Z))\|_{\Theta_V} \leq C \|V\|_{\Theta_V}\|Z\|_{\Theta_Z}.
$$
\end{proposition}

\begin{proof}

The proof relies on a duality argument. We first prove the $L^2B_2^{-1/2,0}$ estimate. Let $U \in L^2,B_2^{1/2,0}$. We write $U = U_1 + U_2$ with 
$$
U_1 = \sum_{j<0} U_j \textrm{ and } U_2 = \sum_{j\geq 0} U_j
$$
where we use the Littlewood-Paley decomposition of $U$. We have, integrating by parts:
$$
\an{U_1,\E(\bar Y W_V(Z)) } = \E\Big( \an{\int_{s}^\infty S (s - t)U_1(t)Y(t) dt, (w*V) Z}\Big).
$$
This yields
$$
\an{U_1,\E(\bar Y W_V(Z)) } \leq \big\|\int_{s}^\infty S_\eta(s - t)U_1(t)Y(t) dt\big\|_{L^{q_1}_{t,x},L^2_\omega}\| (w*V) Z \|_{L^{q_1'}_{t,x},L^2_\omega}.
$$
with $q_1 = d+2$. We apply Lemma \ref{lem:l3a} with $\sigma_1=s = \frac{d}{2}-1$:
$$
\big\|\int_{s}^\infty S (s - t)U_1(t)Y(t) dt\big\|_{L^{q_1}_{t,x},L^2_\omega}\leq \|U_1\|_{L^2_t,B_2^{\sigma_1-\frac12, \sigma_1}}
$$
and since $U_1$ contains only the low frequencies and $s_1 - \frac12 \geq \frac12$, we get
$$
\big\|\int_{s}^\infty S (s - t)U_1(t)Y(t) dt\big\|_{L^{q_1}_{t,x},L^2_\omega}\leq \|U\|_{L^2_t,B_2^{\frac12, 0}}.
$$
By H\"older inequality, since $q_1=d+2$, we get
$$
\| (w*V) Z\|_{L^{q_1'}_{t,x},L^2_\omega}\leq \|w*V\|_{L^2_t,L^2_x}\|Z\|_{L^{p}_t,L^{p}_x,L^2_\omega}.
$$
The estimates above imply for the low frequency part:
\be \lab{bd:quadinter1}
\left| \an{U_1, \E(\bar Y W_V(Z))} \right| \lesssim \|U\|_{B_2^{1/2,0}} \|V\|_{\Theta_V} \|Z\|_{\Theta_Z}.
\ee
We now deal with high frequencies. We have 
$$
\left| \an{U_2,\E(\bar Y W_V(Z)) } \right| = \E\Big( \an{\int_{s}^\infty S(s - t)U_2(t)Y(t) dt, (w*V) Z}\Big).
$$
We get by H\"older inequality
$$
\left| \an{U_2,\E(\bar Y W_V(Z)) } \right| \leq  \big\|\int_{s}^\infty S(s - t)U_2(t)Y(t) dt\big\|_{L^p_t,L^p_x,L^2_\omega} \|V\|_{L^{(d+2)/2}_t,L^{(d+2)/2}_x} \|Z\|_{L^p_t,L^p_x,L^2_\omega}.
$$
From the above estimate we get the control for the high frequency part, by applying lemma \ref{lem:l3a} with $\sigma_1=0$,
\be \lab{bd:quadinter2}
\left| \an{U_2,\E(\bar Y W_V(Z)) }\right| \lesssim \|U_2\|_{B_2^{-1/2,0}} \|V\|_{\Theta_V} \|Z\|_{L^p_t,L^p_x,L^2_\omega}\lesssim  \|U\|_{B_2^{1/2,0}} \|V\|_{\Theta_V}\|Z\|_{L^p_t,L^p_x,L^2_\omega}
\ee
where we used the fact that $U_2$ contains only high frequencies. Combining \fref{bd:quadinter1} and \fref{bd:quadinter2} gives
\be \lab{bd:quadinter3}
\|\E(\bar Y W_V(Z))\|_{L^2_t,B_2^{-1/2,0}} \leq C(f,w) \|V\|_{\Theta_V} \|Z\|_{\Theta_Z}.
\ee
It remains to prove that $\E(\bar Y W_V (Z))$ belongs to $L^{(d+2)/2}_t,L^{(d+2)/2}_x$. Note that $\E(\bar Y W_V(Z))$ belongs to $L^{p_2}_t,L^{q_2}_x$ if $p_2>2,q_2\geq 2$ and
$$
\frac{d}2 - s_2 := \frac2{p_2} + \frac{d}{q_2}  \in [\frac{d}2 - s,\frac{d}2].
$$
This is due to the fact that $Y \in L^\infty_x,L^2_\omega$ and by Strichartz inequalities
$$
\|W_V(Z)\|_{L^{p_2}_t,L^{q_2}_x,L^2_\omega} \lesssim \|V\|_{L^{(d+2)/2}_{t,x}} \|Z\|_{L^p_t,W^{s_2,p}, L^2_\omega}\lesssim \|V\|_{\Theta_V} \|Z\|_{\Theta_Z}
$$
and that $s_2 \in [0,s]$. We recall that by definition $\frac{d}2 -s = 1$. We have
$$
\frac4{d+2} + \frac{2d}{d+2} = 2 \in [1,\frac{d}2],
$$
thus $\E(\bar Y W_V (Z))$ belongs to $L^{(d+2)/2}_t,L^{(d+2)/2}_x$ with:
\be \lab{bd:quadinter4}
\| \E(\bar Y W_V (Z)) \|_{L^{\frac{d+2}{2}}_{t,x}}\lesssim \|V\|_{\Theta_V} \|Z\|_{\Theta_Z}.
\ee
Gathering \fref{bd:quadinter3} and \fref{bd:quadinter4} one obtains the desired continuity estimate
$$
\|\E(\bar Y W_V(Z))\|_{\Theta_V} \leq C(f,w) \|V\|_{\Theta_V}\|Z\|_{\Theta_Z}.
$$

\end{proof}


\section{A space for the initial datum} \lab{sec:ini}

One has the following compatibility result between the space for the perturbation at initial time, and the leading order term for the solution and the potential as given in \fref{id:cz0}. Here, we prove assumption 2 in Proposition \ref{prop-assumption} for dimension higher than $4$ but the proof can be adapted to dimension $3$.

\begin{lemma}

There exists a universal constant $C>0$ such that for all $Z_0\in L^2_\omega ,H^{d/2-1}\cap L^{2d/(d+2)}_x,L^2_\omega $, one has $C_{Z_0}\in \Theta_Z\times \Theta_V$ with
\be \lab{bd:cz0}
\| C_{Z_0} \|_{\Theta_Z\times \Theta_V} \leq C \left( \|Z_0\|_{L^2_\omega ,H^{\frac d2-1}}+\|Z_0\|_{ L^{\frac{2d}{d+2}}_x, L^2_\omega}\right).
\ee

\end{lemma}

\begin{proof} We have that $S(t)Z_0$ belongs to $\Theta_Z$ because of Strichartz estimates and that $Z_0 \in L^2_\omega, H^s$.

Recall the definition of $\Theta_V$, \fref{def:thetav}. The control of the space-time Lebesgue norm of $2\Re \E(\bar Y S(t) Z_0)$ uses standard Strichartz estimates, while the control on its Besov-type norm involves some extra dispersion in the interaction with $Y$. 

We start with the  space-time Lebesgue norm. Since $Y$ belongs to $L^\infty_x,L^2_\omega$ we get by Cauchy-Schwarz:
$$
\| \E(\bar Y S(t) Z_0) \|_{L^{\frac{d+2}{2}}_{t,x}} \lesssim \| S(t)Z_0 \|_{L^{\frac{d+2}{2}}_{t,x},L^2_\omega}.
$$
We recall that 
$$
\frac{4}{d+2} + \frac{2d}{d+2} = 2 \in [1,\frac{d}{2}],
$$
hence by Strichartz estimates, we obtain the first continuity estimate:
\begin{equation}\label{bd:cz0V1}
\| \E(\bar Y S(t) Z_0) \|_{L^{\frac{d+2}{2}}_{t,x}}  \lesssim \| Z_0 \|_{L^2_\omega,H^s_{x}}.
\end{equation}

For the Besov norm, by duality one has:
$$
\| \E(\bar Y S(t) Z_0) \|_{L^{2}_{t},B_2^{-\frac 12,0}} = \sup_{1=\|U\|_{L^2_t ,B_2^{1/2,0}}} \left| \int_{t,x,\omega} S(-t)(\bar Y U)Z_0 \right|.
$$
We write $U=U_1+U_2$ where $U_1=P_{|\xi|\leq 1}U$. For the low frequency part we apply Lemma \ref{lem:l3a} with $p_1=+\infty$, $q_1=2d/(d-2)$ and $\sigma_1=1$:
$$
\| \int_{t} S(-t)(\bar Y U_1) \|_{L^{\frac{2d}{d-2}}_x,L^2_\omega} \lesssim \| U_1 \|_{L^2_t,B_2^{\frac 12,1}}\lesssim \| U \|_{L^2_t,B_2^{\frac 12,0}}.
$$
For the high frequency part, we apply Lemma \ref{lem:l3a} with $p_1=+\infty$, $q_1=2$ and $\sigma_1=0$:
$$
\| \int_{t} S(-t)(\bar Y U_2) \|_{L^{2}_x,L^2_\omega} \lesssim \| U_2 \|_{L^2_t,B_2^{-\frac 12,0}}\lesssim \| U \|_{L^2_t,B_2^{\frac 12,0}}.
$$
Therefore, one has that $\int_t S(-t)(\bar Y U)\in L^{2d/(d-2)}_x,L^2_\omega+L^2_x,L^2_\omega$ (endowed with the canonical norm for sums of Banach spaces) with:
$$
\| \int_t S(-t)(\bar Y U)\|_{L^{2d/(d-2)}_x,L^2_\omega+L^2_x,L^2_\omega} \lesssim \| U \|_{L^2_t,B_2^{\frac 12,0}}.
$$
Hence, by duality:
$$
\left| \int_{t,x,\omega} S(-t)(\bar Y U)Z_0 \right| \lesssim  \| U \|_{L^2_t,B_2^{\frac 12,0}} \left(\|Z_0 \|_{L^{\frac{2d}{d+2}}_x,L^2_\omega}+\|Z_0 \|_{L^{2}_x,L^2_\omega}  \right)
$$
Therefore, by duality one obtains the second estimate:
\be \lab{bd:cz0V2}
\| \E(\bar Y S(t) Z_0) \|_{L^{2}_{t},B_2^{-\frac 12,0}} \lesssim \|Z_0 \|_{L^{\frac{2d}{d+2}}_x,L^2_\omega}+\|Z_0 \|_{L^{2}_x,L^2_\omega} .
\ee
The bounds \fref{bd:cz0V1} and \fref{bd:cz0V2} then imply the continuity estimate:
$$
\| \E(\bar Y S(t) Z_0) \|_{\Theta_V} \lesssim \|Z_0\|_{H^{\frac d2-1},L^2_\omega}+\|Z_0 \|_{L^{\frac{2d}{d+2}}_x,L^2_\omega}.
$$
The identity \fref{id:cz0}, the above bound and \fref{bd:cz01} yield the desired result \fref{bd:cz0}.

\end{proof}


\section{Proof of Theorem \ref{thmmain}}\label{sec-proofth}

We recall that the proof of Theorem \ref{thmmain} relies on finding a solution to the fixed point equation \fref{eq:pointfixe} for the perturbation $Z$ and the induced potential $V$. According to \fref{def:az0}, the fixed point equation can be written in the form:
$$
\mat{Z}{V} = \left(\text{Id}-L\right)^{-1}\left(C_{Z_0} + Q\mat{Z}{V}\right).
$$
This is now a standard routine to solve the above equation thanks to the various estimates derived previously. To solve the above equation, one defines $\Phi[Z_0](Z,V)$ as the mapping
$$
\ba{l l l l}
\Phi[Z_0] : & \Theta_Z \times \Theta_V & \rightarrow &  \Theta_Z \times \Theta_V  \\
 & \mat{Z}{V}  & \mapsto & \left(\text{Id}-L\right)^{-1}\left(C_{Z_0} + Q\mat{Z}{V}\right).
\ea
$$
Let us denote by $\Theta_0=L^2_\omega H^{d/2-1}\cap L^{2d/(d+2)}_xL^2_\omega$ the space for the initial datum with associated norm $\| \cdot \|_{\Theta_0}$. We claim that for $Z_0$ small enough, the mapping $\Phi[Z_0]$ is a contraction on $B(0,C\| Z_0\|_{\Theta_0})$. Indeed, the identity \fref{id:1-l} and the continuity results of Propositions \ref{prop-cont}, \ref{pr:cont1-L1} and \ref{pr:bd 1-L1 L2} give that:
$$
\left(\text{Id}-L\right)^{-1}=\begin{pmatrix} 1 & L_2(1-L_1)^{-1} \\ 0 & (1-L_1)^{-1} \end{pmatrix} \in \mathcal L (\Theta_Z\times \Theta_V).
$$
Hence, for the leading order part, from this and the bound \fref{bd:cz0} for $C_{Z_0}$ one obtains:
$$
\|\left(\text{Id}-L\right)^{-1}\left(C_{Z_0}\right) \|_{\Theta_Z \times \Theta_V}\lesssim \| Z_0\|_{\Theta_0}.
$$
For the quadratic part, recall \fref{id:Q}. In particular, one has
$$
Q(Z,V) = \begin{pmatrix} 0 \\ \E(|Z|^2) \end{pmatrix} + Q_1(Z,V)
$$
where $Q_1$ is bilinear. From Propositions \ref{pr:contquadraZ}, \ref{pr:contWVZ} and \ref{pr:contEWVZ}, if $(Z,V)\in B(0,C\| Z_0\|_{\Theta_0})$ then:
$$
\| Q(Z,V) \|_{\Theta_Z\times \Theta_V} \leq C\|(Z,V)\|_{\Theta_Z\times \Theta_V}^2  \leq  C \| Z_0\|^2_{\Theta_0}
$$
and, due to bilinearity:
\bee
\| Q(Z,V)-Q(Z',V') \|_{\Theta_Z\times \Theta_V} & \leq & C\Big(\|(Z,V)\|_{\Theta_Z\times \Theta_V}+\|(Z',V')\|_{\Theta_Z\times \Theta_V}\Big)\|(Z-Z',V-V')\|_{\Theta_Z\times \Theta_V} \\
&\leq & C \| Z_0\|_{\Theta_0}\|(Z-Z',V-V')\|_{\Theta_Z\times \Theta_V}.
\eee
From the above estimates, on gets that $\Phi[Z_0]$ is indeed a contraction on $B(0,C\| Z_0\|_{\Theta_0})$ for some universal $C$ if $\| Z_0\|_{\Theta_0}$ is small enough. Applying Banach's fixed point theorem yields the existence and uniqueness of a solution to \fref{eq:pointfixe} in $B(0,C\| Z_0\|_{\Theta_0})$. To prove the scattering result, one rewrites \fref{id:Zt} as:
\bee
Z(t) & = & S(t)\left(Z_0 -i\int_{0}^{+\infty} S(-s) (w*V(s)) Z(s) ds -i\int_{0}^{+\infty} S(-s) (w*V(s)) Y(s) ds\right)\\
&&+i\int_t^{+\infty} S(t-s) (w*V(s)) Z(s) ds +i\int_{t}^{+\infty} S(-s) (w*V(s)) Y(s) ds .
\eee
Applying Proposition \ref{pr:contWVZ}, one obtains that $\int_{0}^{+\infty} S(-s) (w*V(s)) Z(s) ds\in L^2_\omega H^{d/2-1}$ and that:
$$
\| \int_t^{+\infty} S(t-s) (w*V(s)) Z(s) ds \|_{L^2_\omega H^{d/2-1}} \lesssim \| Z(s)1_{s \geq t} \|_{\Theta_Z}\| V(s)1_{s \geq t} \|_{\Theta_V}\rightarrow 0
$$
as $t\rightarrow +\infty$. Similarly, from Proposition \ref{prop-cont} one has that $\int_{0}^{+\infty} S(-s) (w*V(s)) Y(s) ds \in H^{d/2-1}$ and that 
$$
\| \int_{t}^{+\infty} S(-s) (w*V(s)) Y(s) ds\|_{L^2_\omega H^{d/2-1}}\lesssim \|_{\Theta_Z}\| V(s)1_{s \geq t} \|_{\Theta_V}\rightarrow 0
$$
as $t\rightarrow +\infty$. Therefore, on has indeed that there exists $Z_{\infty}\in H^{d/2-1}L^2_\omega$ such that, as $t\rightarrow +\infty$:
$$
Z(t)=S(t)Z_{+\infty}+o_{H^{d/2-1}L^2_\omega}(1).
$$
This ends the proof of Theorem \ref{thmmain}.

\section{Example of instability for a rough partition function} \lab{sec:2w}

We study here the linearisation of the dynamics near the superposition of two waves which are orthogonal in probability and propagate in opposite directions. This corresponds to an equilibrium of the form \fref{id:Y} with a rough partition function $f$. Even in the defocusing case, the form of $f$ is involved to ensure linear stability. Indeed, we will prove linear instability for the present example.\\

\noindent Consider a potential $w$ satisfying, without loss of generality if the equation is defocusing, $\int_{\mathbb R^d} w=1$, a mass $m\geq 0$ and a frequency $\xi \in \mathbb R^d$. Let two functions in probability $g_i:\Omega \rightarrow \mathbb C$ for $i=1,2$, with $\int_\Omega | g_i |^2(\omega)d\omega=1/2$ and $\int_{\Omega} \bar g_1(\omega) g_2(\omega)d\omega=0$. The following function is a solution of \fref{eqonrv}:
$$
Y[m,k](\omega,t,x):=\sqrt m e^{-i(|\xi|^2+m)t} \left( g_1(\omega)e^{i \xi.x}+g_2(\omega) e^{-i \xi.x}\right),
$$
which is not stationary, but is at equilibrium. We study a perturbation under the form $X=Y+Z$ and decompose:
$$
Z(\omega,t,x):=g_1(\omega)e^{i(\xi.x-(|\xi|^2+m)t)}\e_1(x,t)+g_2(\omega) e^{i(-\xi.x-(|\xi^2|+m)t)} \e_2 (x,t)+\e_3(\omega,t,x),
$$
where for almost all $t,x\in I \times \mathbb R^d$, one has $\int_{\Omega} \e_3(\omega,t,x)g_i(\omega)d\omega=0$ for $i=1,2$. At the linear level, $\e_1$ and $\e_2$ do not interact with $\e_3$ and their evolution equation is:
\bee
\begin{pmatrix} \pa_t \e_1 \\ \pa_t \e_2 \end{pmatrix} &=& \begin{pmatrix} i\Delta \e_1-2\xi.\nabla \e_1-im\left(w*\Re (\e_1)+w*\Re (\e_2) \right) \\ i\Delta \e_2+2\xi .\nabla \e_2-im\left(w*\Re (\e_1)+w*\Re (\e_2) \right)  \end{pmatrix} +\begin{pmatrix}-ie^{-i(\xi.x-(|\xi|^2+m)t)}\E \left(\bar g_1NL \right) \\ -i e^{(\xi.x+(|\xi|^2+m)t)}\E \left(\bar g_2NL \right)  \end{pmatrix} ,
\eee
where the nonlinear term is $NL:=2(w*\Re (\E (\bar Y Z) ))Z+(w* \E (|Z|^2))(Y+Z)$. We now focus on the linearised operator for $(\e_1,\e_2)$. We decompose between real and imaginary parts, writing $u_1=\Re \e_1$, $u_2=\textrm{Im } \e_1$, $u_3=\Re \e_2$ and $u_2=\textrm{Im } \e_2$. One has the identity
$$
\begin{pmatrix} i\Delta \e_1-2\xi .\nabla \e_1-im \left(w*\Re (\e_1)+w*\Re (\e_2) \right) \\ i\Delta \e_2+2\xi .\nabla \e_2-im\left(w*\Re (\e_1)+w*\Re (\e_2) \right) \end{pmatrix} = \begin{pmatrix} -\Delta u_2-2\xi.\nabla u_1+i(\Delta u_1-2\xi .\nabla u_2-m w*u_1-m w*u_3) \\ -\Delta u_4+2\xi .\nabla  u_3+i(\Delta u_3+2\xi . \nabla u_4-mw*u_1-mw*u_3) \end{pmatrix}.
$$
Consequently the linear coupled dynamics for $\e_1$ and $\e_2$ can be written as the following system:
\be \lab{eq:def A}
\pa_t \begin{pmatrix} u_1 \\ u_2 \\ u_3 \\ u_4 \end{pmatrix} = A   \begin{pmatrix} u_1 \\ u_2 \\ u_3 \\ u_4 \end{pmatrix} , \ \ A:=  \begin{pmatrix} -2\xi .\nabla  & -\Delta & 0 & 0 \\ \Delta -m w* & -2\xi .\nabla & -mw* & 0 \\ 0 & 0 & 2\xi .\nabla & -\Delta \\ -m w* & 0 & \Delta -mw* & 2\xi.\nabla \end{pmatrix}.
\ee
Note that $A$ is a matrix of Fourier multipliers. We now study its spectrum. In the particular case $m=0$, we retrieve for $A$ the block diagonal form
$$
 A:=  \begin{pmatrix} -2\xi .\nabla  & -\Delta & 0 & 0 \\ \Delta & -2\xi .\nabla & 0 & 0 \\ 0 & 0 & 2\xi .\nabla & -\Delta \\  0 & 0 & \Delta & 2\xi.\nabla \end{pmatrix}.
$$
Each of the two matrix operators only have imaginary spectrum, and correspond to a linear Schr\"odinger equation in a moving frame. In the particular case $\xi=0$, $A$ and its symbol are given by:
$$
A:=  \begin{pmatrix} 0  & -\Delta & 0 & 0 \\ \Delta -m w* & 0 & -mw* & 0 \\ 0 & 0 & 0 & -\Delta \\ -m w* & 0 & \Delta -mw* & 0 \end{pmatrix}, \ \ m_A(k)=  \begin{pmatrix} 0  & |k|^2 & 0 & 0 \\ -|k|^2 -m \hat w(k) & 0 & -m\hat w(k) & 0 \\ 0 & 0 & 0 & |k|^2 \\ -m \hat w(k) & 0 & -|k|^2 -m\hat w(k) & 0 \end{pmatrix}.
$$
The eigenvalues of $m_A$ are by a direct check $\lambda_{\pm,\pm} =\pm \sqrt{-|k|^4-(1\pm 1)m|k|^2 \hat w (k)}\in i \mathbb R$ in the defocusing case $\hat w\geq 0$. This linear operator is similar to the one arising in the linearisation of the Gross-Pitaevskii equation near the trivial state, and has thus in the defocusing case only imaginary spectrum. The situation is more involved when $\xi\neq 0$ and $m\neq 0$. In the case of a Dirac potential one has the following instability result, whose proof can be extended to include other additional potentials.

\begin{lemma}

Let $\hat w=1$, $m>0$ and $\xi \neq 0$. Then the Fourier multiplier of the differential matrix $A$ given by \fref{eq:def A} possesses positive and negative eigenvalues in the vicinity of the frequency $k=\xi \sqrt{4-\min(2,m/|\xi|^2) } $.

\end{lemma}

\begin{proof}

We compute the characteristic polynomial of $A$, using the notations $a=2\xi .\nabla$, $b=-\Delta$ and $c=mw*$ to ease computations:
\bee
P_A(\mathcal X) &:=& \det (\mathcal X\text{Id}-A)= \det \left( \begin{pmatrix} \mathcal X+2\xi .\nabla  & \Delta & 0 & 0 \\ -\Delta +mw* & \mathcal X+2\xi .\nabla  & mw* & 0 \\ 0 & 0 & \mathcal X-2\xi .\nabla & \Delta \\ mw* & 0 & -\Delta +mw* & \mathcal X-2\xi .\nabla \end{pmatrix} \right)\\
&=& \mathcal X^4+2((b+c)b-a^2)\mathcal X^2+((b+c)b+a^2)^2-b^2c^2.\\
\eee
We set $\mathcal Y=\mathcal X^2$ and compute the discriminant:
$$
4D^2=(2(b+c)b-a^2)^2-4\left(((b+c)b+a^2)^2-b^2c^2\right)=4b\left(bc^2-4(b+c)a^2 \right).
$$
Therefore the roots in $\mathcal Y$ of the above polynomial are:
$$
\mathcal Y_{\pm}:= a^2-(b+c)b\pm D=4(\xi .\nabla )^2+(-\Delta+m w*)\Delta\pm \sqrt{-\Delta \left(-\Delta (mw*)^2-16(-\Delta+mw*)(\xi .\nabla )^2 \right)}
$$
and the roots in $\mathcal X$ are:
$$
\mathcal X_{\pm,\pm}=\pm \sqrt{\mathcal Y_{\pm}}.
$$
If $\mathcal Y_{+}$ or $\mathcal Y_-$ has some positive spectrum, then at least one of the $\mathcal X_{\pm,\pm}$ has some positive and negative spectrum which signals a linear instability. For $r>0$, the symbol associated to $\mathcal Y_{+}$ evaluated at $r\xi$ is:
\bee
&& -4r^2|\xi|^4-r^4|\xi|^4-m |\xi|^2r^2 +\sqrt{ r^2|\xi|^2 \left(r^2|\xi|^2m^2+16(r^2|\xi|^2+m)|\xi|^4r^2\right)}\\
&=& |\xi|^2r^2\left( -4|\xi|^2-r^2|\xi|^2-m +\sqrt{m^2+16r^2|\xi|^4+16m|\xi|^2}  \right).
\eee
The above quantity is positive when the following polynomial
$$
\left(4|\xi|^2+r^2|\xi|^2+m\right) -\left(m^2+16r^2|\xi|^4+16m|\xi|^2 \right)= |\xi|^4(r^2-4)\left(r^2-4+2\frac{m}{|\xi|^2}\right)
$$
is negative. Therefore, at the frequency $k=\xi \sqrt{4-\min(2,m/|\xi|^2) } $, the Fourier multiplier of $A$ has a positive and a negative eigenvalue.

\end{proof}

\end{document}